\documentclass[11pt]{amsart}
\usepackage[english]{babel}

\usepackage{amsmath,amsthm,amsfonts,amssymb}
\usepackage[foot]{amsaddr}
\usepackage{mathtools}
\usepackage{color}

\usepackage{booktabs} 
\usepackage{amssymb}
\usepackage{amsfonts}
\usepackage{amsmath}
\usepackage{mathrsfs}
\usepackage{mathbbol}
\usepackage{amsfonts}
\usepackage{relsize}
\usepackage{booktabs}
\usepackage{stackrel}
\usepackage{mathrsfs}
\usepackage{mathbbol}

\usepackage{mathbbol}

\usepackage[active]{srcltx}
\usepackage{enumitem}

\def\text{\mbox}

\def\div{\mathrm{div}\,}

\def\dps{\displaystyle}
\def\R{{\mathbb R}}

\def\CC{{\mathscr C}}

\def\pt{\partial}

\def\<{\langle}
\def\>{\rangle}

\def\ker{N}
\def\bve{\;|\;}

\newenvironment{remark}{ {\sc Remark -- } }   {\\}  
\newtheorem{theorem}{Theorem}[section]
\newtheorem{proposition}[theorem]{Proposition}
\newtheorem{lemma}[theorem]{Lemma}
\newtheorem{corollary}[theorem]{Corollary}

\newcommand{\bfgreek}[1]{\bm{\@nameuse {up#1}}}

\newcommand{\vc}[1]{ {#1}}

\def\Ã«{\vc{\xi}}

\def\hf{{\widehat{f}}}

\newcommand{\vecc}[1]{{#1}}

\def\u{\vecc{u}}

\def\v{\vecc{v}}
\def\w{\vecc{w}}

\def\zero{\vecc{0}}

\def\e{\vecc{e}}

\def\x{\vecc{x}}
\def\y{\vecc{y}}

\newenvironment{ourproof}[1]{ \begin{quote}{\sc Proof  #1} -- }   {  $\blacksquare$  \end{quote}   }


\def\om{\Omega}
\def\dps{\displaystyle}

\title{}

\def\z{z}
\def\x{x}
\def\y{y}
\def\v{v}

\def\RDC{{C}}

\def\argmin{{\rm argmin}\,}

\renewcommand{\leq}{\leqslant}
\renewcommand{\geq}{\geqslant}

\newcommand{\drv}[1]{\dot{#1}}
\newcommand{\ddrv}[1]{\ddot{#1}}
\newcommand{\card}[1]{{|#1|}} 
 
\newcommand{\TG}[2]{T(#1, #2)}
\newcommand{\NC}[2]{N(#1, #2)}

\def\ggrec{\boldsymbol}
\def\teps{{\ggrec \varepsilon}}
\def\tepsE{\teps^{\rm e}}
\def\tepsP{\teps^{\rm p}}
\def\drvepsE{\drv{\teps}^{\rm e}}
\def\drvepsP{\drv{\teps}^{\rm p}}
\def\tsig{{\ggrec \sigma}}

\def\tsigs{{\ggrec  \sigma}^\star}
\def\ttau{{\ggrec \tau}}
\def\tita{{\ggrec  \eta}}
\def\tkappa{{\ggrec  \kappa}}
\def\tvarphi{{\ggrec  \varphi}}
\def\tstrr{{\ggrec \varepsilon}}
\def\strr{{\varepsilon}}

\def\sig{ \sigma}
\def\tsigs{{\ggrec \sigma^\star}}

\def\eps{\varepsilon}

\def\dvsig{{\ggrec {\overline{\sigma}}}}

\def\dvsig{\overline{\tsig}}

\def\dvtau{\overline{\ttau}}
\def\lda{ \lambda}
\def\trc{{\rm tr}}

\def\id{{\rm{\bf{I}}}}
\def\IDT{{\mathbb{Id}}}
\def\tgrad{{\boldsymbol \nabla}}

\def\symt{{\mathbb M}^{3\times 3}_{\rm{sym}}}
\def\symtp{{\mathbb M}^{3\times 3}_{\rm{sym, +}}}
\def\symtn{{\mathbb M}^{n\times n}_{\rm{sym}}}
\def\symtpn{{\mathbb M}^{n\times n}_{\rm{sym, +}}}
\def\nsymt{{\mathbb M}^{3\times 3}}
\def\Dvs{{\mathbb M}^{3\times 3}_{\rm{D}}}
\def\Hyd{{\R \id}}
\def\span{{\rm{span}}}

\def\LDA{\Lambda}

\newcommand{\smatsym}[6]{ 
\left[
\begin{array}{ccc}
\dps{#1 }  & \dps{#4 }  & \dps{#5 }  \\
\dps{#4}  & \dps{#2}  & \dps{#6 }  \\
\dps{#5}  & \dps{#6}  & \dps{#3}
\end{array}
\right]
}

\newcommand{\diagn}[3]{\smatsym{#1}{#2}{#3}{0}{0}{0}}

\def\Iv{I}
\def\IJ{J}

\title{}
\def\HILB{{\mathscr H}}
\def\HOOK{{\mathbb H}}
\def\MHOOK{{\boldsymbol {\mathscr H}}}

\newcommand{\prj}[1]{{\Pi}_{#1}}
\newcommand{\PRJ}[1]{\Pi_{#1}}

\newcommand{\parallelsum}{\small{\mathbin{\!/\mkern-5mu/\!}}}
\newcommand{\Sat}[1]{{\rm{Sat}}(#1)}

\def\ES{{\mathscr N}}
\newcommand{\PRJP}[1]{\Pi_{#1}^{\mathscr N}}
\newcommand{\PRJO}[1]{\Pi_{#1}^{T}}

\def\thet{\theta}

\def\ker{N}

\def\ftrs{f_{{T}}}

\def\CTRS{C_{T}}
\newcommand{\ORT}[1]{{\bold O}_{#1}}
\def\SSD{{\mathbb S}}
\def\force{h}
\def\KS{K}

\newcommand{\Diag}[1]{{\rm Diag}(#1)}

\title[Model of elastic perfectly plastic solids:  a convex analysis approach]{Unification of the mathematical model of elastic perfectly plastic solids:  a convex analysis approach}

\author[T. Z. Boulmezaoud]{Tahar Z. Boulmezaoud$^{1, 2, 3}$}
\address{\rm  $^1$ Universit\'e Paris-Saclay, UVSQ, LMV, Versailles, France.}
\address{\rm  $^2$ Centre National de la Recherche Scientifique (CNRS), D\'el\'egation Paris Michel-Ange, 3 rue Michel-Ange, 75016 Paris, France.}
\email{tahar.boulmezaoud@uvsq.fr}
\author[B. Khouider]{ Boualem Khouider$^{3}$}
\address{\rm  $^3$ Mathematics and Statistics, University of Victoria, PO BOX 1700 STN CSC, Victoria, B.C.,  Canada V8W 2Y2. }
\email{khouider@uvic.ca}


 
\begin{document}

\maketitle

\begin{abstract}
A new mathematical formulation for the constitutive laws governing elastic perfectly plastic materials is proposed here. In particular, it is shown that the elastic strain rate and the plastic strain rate form an orthogonal decomposition with respect to the
tangent cone and the normal cone of the yield domain. It is also shown that the stress rate can be seen as the projection on the tangent cone of the elastic stress tensor. This approach leads to a coherent mathematical formulation of the elasto-plastic laws
 and simplifies the resulting system for the associated flow evolution equations. The cases of one or two yields functions are treated in detail. 
 The practical examples of the von Mises and  Tresca yield criteria are worked out in detail to demonstrate the usefulness of the new formalism  in applications. 
 \end{abstract}

\tableofcontents

\section{Introduction}
In the theory of perfect plasticity, the deformation of a material is mainly decomposed into two components; 
an elastic deformation due to reversible microscopic processes for which there is a one-to-one map between 
the stress and the strain, and an irreversible plastic deformation  for which such a map is lost.
 The calculation of deformations of an elasto-plastic material must therefore take into account not only
its current state, but also its history. For this reason, the common approach consists to  find this 
deformation as the result of infinitesimal variations of the strain and the stress tensors  (see, e. g., 
\cite{lubliner1970},  \cite{hashi2012}, \cite{hanreddy}).  \\
$\;$\\
Before describing the aim and main results of this work, let us start with a simple reminder of the general principles of elasto-plasticity.   Leaving aside thermal effects and hardening  (the latter will be the subject of a specific 
study in a future paper), the deformation of a material occupying - in its undeformed state -
a domain $\om \subset \R^3$ is described by the knowledge of 
the displacement vector field $\u(\x, t)$ characterizing the geometry,
 and the Cauchy stress tensor $\tsig(\x, t)$ characterizing  the state of
the material (with $x$ is any point of $\om$ and $t$ is time). In incremental 
elastic perfectly plastic model, the displacement $\u$ (assumed to be small) and the stress
 $\tsig$ are governed by the following usual 
principles:
\begin{itemize}
\item {\it (additive decomposition)}  The strain rate tensor $\drv{\tstrr}$  is the sum the elastic strain rate $\drv{\tepsE} $ and the plastic strain rate $\drv{\tepsP}$: 
\begin{equation}\label{sumdecompo0}
\drv{\tstrr}= \drv{ \tepsE} +\drv{\tepsP}.
\end{equation}
\item {\it (the yield criterion) } The stress tensor satisfies 
\begin{equation}\label{yield_principle}
\tsig \in \RDC, 
\end{equation}
 where $ \RDC $ is a nonempty closed {\it convex} subset of three-by-three symmetric tensors, $\symt$ (see, e. g., \cite{drucker1949}, \cite{hill56}, \cite{lubliner1970}).  
It is assumed that the material is perfectly plastic, that is $\RDC$ is constant during the deformation process (no hardening or softening occurs).  When the stress is strictly inside $\RDC$,  the strain-rate is equal 
to the elastic stress-rate $\drv{\tstrr} = \drv{\tepsE}$ and  $\drv{\tepsP} = 0$. The plastic onset 
occurs when the stress reaches the boundary  $\pt \RDC$ of $\RDC$ ($\pt \RDC$ is called 
the yield surface). Note that many yield criteria used in practice are commonly defined by functional constraints of the form
\begin{equation}\label{inequal_yield}
f_i(\tsig) \leq 0 \mbox{ for } 1 \leq  i  \leq m, 
\end{equation}
where the yield functions $f_i$ depend only on the principal stresses of $\tsig$. \\ 
There is a large number of criteria describing the yield of materials in the literature. 
 The most commonly used are the Tresca criterion (\cite{tresca1964}) and the von Mises criterion \cite{mises1928}.

\item {\it (Hooke's Law)} The elastic strain rate $\drv{\tepsE}$ is related to the stress rate by 
\begin{equation}\label{hookzero}
\drv{ \tsig} = \HOOK(\drv{\tepsE}), 
\end{equation}
where $\HOOK$ is the fourth-order isotropic elasticity tensor given by  \eqref{law_hooke4}.
\item  {\it (Principle of maximum work)} when $\tsig \in\pt \RDC$,  the pair $(\tsig,  \drvepsP)$ satisfies 
\begin{equation}\label{intro_normal} 
 \drvepsP:\tsig \geq  \drvepsP:{\tsigs} \mbox{ for all } \tsigs \in \RDC. 
\end{equation}
(see  \cite{hill48, hill50},  \cite{koiter1953, koiter1960}, \cite{lubliner1970}). 
\item {\it (Consistency condition) }  when $\tsig \in\pt \RDC$,  $ \drvepsP$ and $ \drv{\tsig}$ verify
\begin{equation}\label{consistancy0} 
 \drvepsP : \drv{\tsig}  = 0,
 \end{equation} 
 at all times.
 
As indicated in \cite{DuvautLions}, condition \eqref{consistancy0} results from 
   \eqref{intro_normal}  under strong time regularity assumptions on $\tsig$ (here 
    $\tsig$ is a tensor valued function depending on time $t$ and position $\x$). 
  This might explain why this condition  
  is often sidelined and not taken into account by many authors. 
  \end{itemize}
In view of \eqref{sumdecompo0} and \eqref{hookzero}, Condition  \eqref{intro_normal} is often written in the following 
form,  called the normality rule,
\begin{equation}\label{flowrule0} 
 \drvepsP \in \NC{\RDC}{\tsig},
\end{equation}
or, equivalently, 
\begin{equation}\label{orthoPDS} 
 \drv{\tstrr}  -  \HOOK^{-1}(\drv{\tsig}) \in \NC{\RDC}{\tsig},
\end{equation}
where $\NC{\RDC}{\tsig}$ denotes the normal cone of $\RDC$ at $\tsig$ and $\HOOK^{-1}$ is
the inverse of the operator  $\HOOK$. \\

Because of the inherent irreversibility of plastic deformations, it is more meaningful 
to describe elasto-plastic deformation processes by their  infinitesimal variations, i.e. by the strain rate and stress rate tensors. 
Besides, these principles are complemented by the local equations governing the displacement of the material elements, given  
in  the general form
\begin{equation}\label{evomotion_eq}
\rho \ddrv{\u} - \div \tsig - \force = 0, 
\end{equation}
where $u$ is the displacement of the material element with respect to its original spacial coordinate $x$ with the two dots on its top denoting the second derivative in time, $ \force$ is the volume density of forces, $\rho$ the density of the material, and $\div \tsig$ the vector field
given by
$
(\div \tsig)_i = \sum_{j=1}^3 \frac{\pt \sigma_{i j}}{\pt x_j} \mbox{ for } 1 \leq i \leq 3,
$
 augmented by appropriate boundary and initial conditions set by the forces acting on the boundary of the material and its initial state. These physical constraint that are nonetheless necessary in order to find a meaningful solution to this system of equations are not important for the present work. The proposed formulation can nevertheless accommodate any type of boundary and initial conditions. 
 
It is worthwhile mentioning that under the assumption of a quasi-static evolution, equation \eqref{evomotion_eq} becomes
\begin{equation}\label{qsmotion_eq}
 \div \tsig + \force = 0, 
\end{equation}
and that this type of assumption is often made in practice when the elasto-plastic times scales of the material are much faster than those of the underlying volume and boundary forces.

The elastic and perfectly plastic time dependent problem of a material occupying a geometrical domain is
  often written as a variational inequality (see, e. g., \cite{DuvautLions}). 
 It leads to a non-linear and cumbersome large number of equations although it has been studied by several authors, both in the case of quasi-static evolution and in the case of full dynamics (see, e. g., \cite{DuvautLions}, \cite{johnson1976},  \cite{temam1986},  \cite{fuchs}, \cite{dalmaso2006}, \cite{baba2021}). 
The present work aims at reformulating the local principles of elasto-plasticity 
into a smaller number of equations,  which allows in particular
 to get rid of the inequalities (and thus of the variational inequalities
   associated with the global time evolution problem).
Our approach is based on the following statement: behind the system lies 
the unique orthogonal decompositions of the strain rate $\drv{\tstrr}$ and its image 
 $\HOOK(\drv{\tstrr})$ in the form $ \ttau+\tita$ with  $\ttau$ belonging
 to the tangent cone of $\RDC$ at $\tsig$ and $\tita$ belonging
 to its polar cone (i. e. the normal cone  of $\RDC$ at $\tsig$). More
 precisely, we shall prove the following.
$$
 \drvepsE =   \PRJ{\TG{\RDC}{\tsig}}   \drv{\tstrr}, \; \;  \drvepsP  = \dps{  \PRJ{\NC{\RDC}{\tsig}}  \drv{\tstrr}}, \; \; \drv{\tsig} =   \PRJ{\TG{\RDC}{\tsig}}\HOOK   (\drv{\tstrr}),
$$
where $ \PRJ{\TG{\RDC}{\tsig}}$ (resp.  $\PRJ{\NC{\RDC}{\tsig}}$) represents the orthogonal 
projection on $\TG{\RDC}{\tsig}$ (resp. on ${\NC{\RDC}{\tsig}}$).  This will allow us in particular to formulate the system as an augmented evolution system of the form
\begin{equation}\label{evol_sys}
\frac{\pt}{\pt t}
\left(
\begin{array}{c}
v \\
\tsig 
\end{array}
\right) -
\left(
\begin{array}{c}
\div \tsig \\
\MHOOK(\tsig, \drv{\tstrr} (v))
\end{array}
\right) = 
\left(
\begin{array}{c}
\force \\
0
\end{array}
\right), 
\end{equation}
with $\MHOOK(\tsig, \drv{\tstrr}(v))  = \PRJ{\TG{\RDC}{\tsig}}\HOOK(\drv{\tstrr}(v))$ and $v=\dot u$ is the flow velocity of the material elements. The nonlinear equality and inequality constraints associated with the yield criterion, which are inevitable when using a variational method for example,  are replaced by easy to calculate projections onto the tangent and normal cones, as demonstrated by the application of this new methodology  to the von Mises and Tresca  criteria. \\
 
Time differentiation  of the first component of the system leads to  an evolution equation involving only the velocity
\begin{equation}
\frac{\pt^2 v }{\pt t^2}  -  \div \MHOOK(\tsig, \drv{\tstrr} (v))  = \frac{\pt h}{\pt t}. 
\end{equation}
From all these elements, it follows that it is necessary to calculate the projectors 
$\PRJ{\TG{\RDC}{.}}$ and  $\PRJ{\NC{\RDC}{.}}$
 on the tangent cone and the normal cone. \\
 This approach is quite general, and does not require assumptions about the regularity of $\RDC$, 
nor restrictions on the number of functions that define it.  However, we will present the 
calculations in the cases where $\RDC$ is defined by one or two yield functions only, to demonstrate the usefulness of the new methodology in practical applications. The case of the Von
 Mises criterion will be  in particular well detailed. The Tresca criterion, for
  which the domain has degenerate corners, will also be carefully examined.  
  Finally, it will be discussed how invariance by similarity of the yield functions can be exploited to compute 
easily and efficiently the projectors involved in the formulation. 

The rest of this paper is organized as follows. This section will end with preliminaries and notations. 
In Section \ref{sec_main_res}, we present the main results allowing to reformulate the equations of the elastic perfectly  plastic model in terms of the projectors on the normal and tangent cones.  In Section \ref{consti_laws}, we will treat in more details the case where we have one or two functions defining the yield surface. Section \ref{VM-T-cr} is devoted to the Von Mises's and Tresca's criteria. 
The last section is devoted to some concluding remarks concerning  invariance by similarity of
yield  functions  and upcoming extensions.

{\it Notations and preliminaries}\label{sect_notat}
In the sequel, all the elements of $\R^3$ will be considered as column vectors.   
For $\x, \, \y \in \R^3$, we denote by $\dps{\x}\cdot{\y} \in \R^3$ their Euclidian scalar product, 
by $\x \otimes \y= \x \y^t$ their tensor product and by  $\x \odot \y= \frac{1}{2}(\x \otimes \y +\y \otimes \x)$
their symmetric tensor product.  The superscript $t$ denotes the matrix or vector transpose.The same notation will be used for any 
second order tensor $\tsig$ and for the matrix of its components $(\sigma_{ij})_{1 \leq i, \; j \leq n}$
 (with respect to the canonical basis of $\R^3$). Given two second order tensors $\tsig$
  and $\teps$, we denote by  $\tsig \teps$ their products, by  $\tsig:\teps = \trc{(\tsig \teps^T)}$ 
  their scalar product and by $\|.\|$ the associated (Frobenius) tensor norm. When   $\tsig$ and $\teps$ are symmetric,
we have 
$$
\tsig : \teps = \sum_{i,j=1}^3 \sigma_{ij} \eps_{ij} = \sum_{i =1}^3 \sigma_{i i}\eps_{ii} + 2 \sum_{1 \leq i < j \leq 3} \sigma_{i j}\eps_{ij},
$$
and 
$$
\|\tsig\| = \{\sigma_{11}^2 +  \sigma_{22}^2 + \sigma_{33}^2 + 2 \sigma_{12}^2 +  2 \sigma_{13}^2 + 2 \sigma_{23}^2 \}^{1/2}. 
$$
The identity tensor of second order (resp. of order four) will be denoted by  $\id$ (resp. $\IDT$). 
Given an integer $n \geq 1$, $\symtn$ denotes the space of $ n\times n$ symmetric matrices,   $\symtpn$ is the subset of $\symtn$ comprised of positive semidefinite matrices and 
 $\ORT{n}$ is the group of $n \times n$ orthogonal matrices. Given a symmetric tensor $\tsig \in \symt$,   $\lda_1(\tsig)$, $\lda_2(\tsig)$ and $\lda_3(\tsig)$ 
denote its principal values (or eigenvalues) 
aranged in decreasing order: $\lda_1(\tsig) \geq \lda_2(\tsig) \geq \lda_3(\tsig)$. We set
$$
\LDA (\tsig) = (\lda_1(\tsig), \lda_2(\tsig), \lda_3(\tsig)).
$$
 We denote by $\Iv_1(\tsig)$, $\Iv_2(\tsig)$ and $\Iv_3(\tsig)$ the principal invariants of $\tsig$ such that
   the characteristic  polynomial of $\tsig$ is given by 
\begin{equation}
 \det(\lambda I - \tsig) = \lambda^3 - \Iv_1(\tsig)\lambda^2 +  \Iv_2(\tsig)\lambda - \Iv_3(\tsig) \mbox{ for all }  \lambda \in \R.
\end{equation}
We have
\begin{equation}
\Iv_1(\tsig) =  \trc(\tsig), \; \Iv_2(\tsig) =  \frac{1}{2} ( \trc(\tsig)^2 - \trc(\tsig^2)), \; \Iv_3(\tsig) =  \det(\tsig).
\end{equation}


If $f\; :\;\symt \to \R$ is differentiable  at $\tsig \in \symt $, then $\tgrad f(\tsig)$  will be its symmetric 
 gradient at $\tsig$,  that is $\tgrad f(\tsig)$  is the unique second order tensor satisfying 
$$
\forall \tsigs \in \symt, \; df(\tsig)\tsigs  = (\tgrad f(\tsig) ) : \tsigs,
$$
where $df(\tsig)$ is the differential of $f$ at $\tsig$.  It is worth noting 
 that when the function $f$ is defined over the nine dimensional space $\nsymt$  of $3 \times 3$
 matrices with real entries,   its symmetric gradient  $\tgrad f(\tsig)$  is different from its 
 gradient as a function on  $\nsymt$ (the former is the symmetric part of the latter). 
 More precisely,  the components of $\tgrad f(\tsig)$  are given by
$$
\tgrad f(\tsig) = \sum_{1 \leq i \leq j \leq 3}   \frac{\pt f}{\pt \sig_{ij}}(\tsig) \e_i  \odot \e_j. 
$$
We introduce the subspace of deviatoric tensors (or matrices): 
\begin{equation}
\Dvs = \{\tsig \in \symt \bve \trc(\tsig) = 0\},  
\end{equation}
Obviously,  $\symt = \Hyd \oplus^\perp \Dvs$ and for all $\tsig \in \symt$, we can write 
\begin{equation}
 \tsig  = \dvsig+ \frac{\trc{(\tsig)}}{3} \id,
\end{equation}
where $\id$ is the identity tensor and $\dvsig \in   \Dvs $ is the deviator of $\tsig$.  Obviously $\Iv_1(\dvsig)=0$ and 
\begin{equation}
\Iv_2(\dvsig) =  \dps{   \Iv_2(\tsig) - \frac{\Iv_1(\tsig)^2}{3}}, \; 
\Iv_3(\dvsig) = \dps{\frac{ \trc(\dvsig^3) }{3} =     \frac{2 \Iv_1(\tsig)^3}{27}  -  \frac{\Iv_1(\tsig) \Iv_2(\tsig)  }{3}+ \Iv_3(\tsig).} 
\end{equation}
It is customary in solid mechanics to set $\IJ_2(\tsig) = - \Iv_2(\dvsig)$ and $\IJ_3(\tsig) = \Iv_3(\dvsig)$. 
We have 
\begin{eqnarray}
\IJ_2(\tsig) & =& \dps{ \frac{1}{6} ((\sig_{11} - \sig_{22})^2 + (\sig_{22} - \sig_{33})^2 + (\sig_{11} - \sig_{33})^2) } \nonumber \\
&& + \sig_{12}^2  +  \sig_{23}^2  +  \sig_{13}^2, \\
 &=& \dps{ \frac{1}{6} ((\lda_{1} - \lda_{2})^2 + (\lda_{2} - \lda_{3})^2  +(\lda_{1} - \lda_{3})^2 ), }  \\
 &=& \frac{1}{2} \|\dvsig\|^2.
\end{eqnarray}
Also, by the Cayley-Hamilton theorem, which states that every matrix is a solution of its characteristic equation, we have 
\begin{equation}
\tsig^3 =  \Iv_1(\tsig)\tsig^2 -  \Iv_2(\tsig)\tsig + \Iv_3(\tsig) \id, 
\end{equation}
and thus \[\dvsig^3 =  \IJ_2(\tsig) \dvsig   +  \IJ_3(\tsig)  \id. \]
In the sequel, we also need to express the gradient of these invariants. The reader can easily verify the
following identities
\begin{eqnarray}
\tgrad \Iv_1(\tsig) &=& \id, \;\\
 \tgrad \Iv_2(\tsig) &=&  \Iv_1(\tsig)  \id- \tsig, \; \\
 \tgrad \Iv_3(\tsig) &=&    \tsig^2 - \Iv_1(\tsig) \tsig + \Iv_2(\tsig) \id, \\
\tgrad \IJ_2(\tsig) &=& \dvsig, \;\\
 \tgrad \IJ_3(\tsig)  &=&  \dvsig^2 -\frac{2}{3}  \IJ_2(\tsig) \id. \label{gradJ2J3}
\end{eqnarray}
Now, given a nonempty convex set $\KS \subset \symt$ and $\tsig\in \KS$, the tangent cone $\TG{\KS}{\tsig}$ to $\KS$ at $\tsig$ is defined by 
\begin{equation}
\TG{\KS}{\tsig} = \overline{\{\alpha (\tita-\tsig)\bve \tita \in \KS \mbox{ and } \alpha > 0\},}\label{tc1}
\end{equation}
where the overline denotes the topological closure of the underlying set. This is a closed convex cone of 
$\symt$.  
The normal cone $\NC{\KS}{\tsig}$ of $\KS$ at $\tsig$ is defined as the dual cone of $\TG{\KS}{\tsig}$, that is
\begin{equation}
\NC{\KS}{\tsig} =  \{\tita \in \symt  \bve  \ttau : \tita \leq 0  \mbox{ for all } \ttau \in \TG{\KS}{\tsig} \}. \label{nc1}
\end{equation}
We recall that $\TG{\KS}{\tsig} = \symt$  and $\NC{\KS}{\tsig} = \{\zero\}$  when $ \tsig$ belongs to the interior of $\KS$. 
We have the two  elementary properties
\begin{itemize}
\item For all $\alpha > 0$, $\TG{\alpha \KS}{\alpha \tsig} = \TG{\KS}{\tsig}$, $\NC{\alpha \KS}{\alpha \tsig} = \NC{\KS}{\tsig}$.
\item For all $\ttau \in \symt$,  $\TG{ \KS+\ttau}{\tsig+\ttau} = \TG{\KS}{\tsig}$,  $\NC{ \KS+\ttau}{\tsig+\ttau} = \NC{\KS}{\tsig}$.
\end{itemize}
Finally, we denote by $\PRJ{\KS}$ the 
orthogonal projection on the convex $\KS$ as an operator on $\symt$.  For $\ttau \in  \symt$, we have
\begin{equation}
 \PRJ{\KS} \ttau = \argmin_{\tita \in \KS }\|\ttau - \tita\|. \label{prc1}
\end{equation}
\section{Reformulating elasto-plasticity equations}\label{sec_main_res}
Consider an elasto-plastic material occupying a region $\om \subset \R^3$.
 When volume and/or surface forces are applied to the material body, the deformation and the state of this material 
can be characterized by the evolution of  the displacement vector field $(x, t) \in \om \times I \mapsto u(x, t)$
 and the Cauchy (internal) stress tensor $(x, t) \in \om \times I \mapsto  \tsig(x, t)$, respectively. Here, $I$ is the time interval during which the deformation take place. For convenience, we assume that I I is semi-open on the form $I=[0,T)$. 
 We denote by $v = \pt u/{\pt t}$ the velocity vector field, and $\drv{\tsig}$ the time derivative of $\tsig$ while $\tstrr$ and $\drv{\tstrr}$ are respectively the strain and the strain rate tensors: 
\begin{equation}\label{strainsrate_def}
  \strr_{i,j} = \frac{1}{2} \left(  \frac{\pt u_i}{\pt x_j}  +   \frac{\pt u_j }{\pt x_i} \right), \;   {\drv{\strr}}_{i,j} = \frac{1}{2} \left(  \frac{\pt v_i}{\pt x_j}  +   \frac{\pt v_j }{\pt x_i} \right), \; 1 \leq i, j \leq 3. 
\end{equation}
 We are interested in the evolution in time of these quantities, locally, at all points
 $\x \in \om$. The focus is on the constitutive law which links 
 the strain, the strain rate, the stress, and the stress rate tensors.  \\

 As stated in the introduction section, for many materials, it is meaningful to assume that the internal stress is restricted to a closed convex set $\RDC$ of $\symt$ and that 
  $\RDC$ is insensitive to hydrostatic pressure, that is 
\begin{equation}\label{condition_indiff_hydro}
\forall \tsig \in \RDC, \forall p \in \R, \; \tsig + p \id \in \RDC.
\end{equation}
Geometrically, this means that $\RDC$ is an infinite cylinder aligned along the identity matrix. In particular,
$\RDC$ is unbounded. The extension to materials not complying with \eqref{condition_indiff_hydro} 
will be the subject of  future work.

This indifference assumption has a direct implication on the normal and the tangent cones. We have
\begin{equation}\label{indif_conesNT}
\Hyd \subset  \TG{ \RDC}{\tsig}, \; \NC{ \RDC}{\tsig} \subset (\Hyd)^\perp
\end{equation}
and, in particular, 
\begin{equation}\label{indif_conesNT2}
\trc(\ttau) = 0 \mbox{ for all }  \ttau \in \NC{ \RDC}{\tsig}. 
\end{equation}
Here and subsequently, we drop the $(x, t)$ arguments for simplicity of notation. \\
 
Inside $\RDC$, the internal stress $\tsig$ is related to the elastic strain through Hooke's law \eqref{hookzero}  with
\begin{equation}\label{law_hooke4}
\HOOK(\ttau) = 2 \mu \ttau+ \lda  \trc(\ttau) \id, \; \; \mbox{ for all } \ttau \in \symt, 
\end{equation}
where $\lambda > 0$ and $\mu >0$ are the Lam\'e coefficients, which are assumed to be constant. We can also write
$\tepsE = \HOOK^{-1}(\tsig)$ with  
\begin{equation}\label{law_hooke1}
 \HOOK^{-1}(\tsig)  = - \frac{\lambda}{2 \mu (3 \lda + 2 \mu)} \trc(\tsig) \id + \frac{1}{2 \mu} \tsig, 
\end{equation}
or
\begin{equation}\label{law_hooke2}
 \HOOK^{-1}(\tsig)  = - \frac{\nu}{E} \trc(\tsig) \id + \frac{1+\nu}{E} \tsig, 
\end{equation}
where $E$  and $\nu$ are respectively the Young's modulus and Poisson's ratio, given by
\begin{equation}\label{Hooke_consts}
E = \frac{ \mu(3 \lda + 2 \mu)}{\lda +\mu}, \;  \nu = \frac{\lda}{2(\lda +\mu)}.
\end{equation}

\begin{theorem}\label{moreau_elastoplast}
 The following two statements are equivalent:
\begin{enumerate}
\item $(\tstrr, \tepsE, \tepsP, \tsig)$ satisfy  \eqref{sumdecompo0},  \eqref{yield_principle},  \eqref {hookzero},  \eqref{intro_normal} and \eqref{consistancy0}.
\item $(\tstrr, \tepsE, \tepsP, \tsig)$ satisfy \eqref{yield_principle} at $t=0$ 
  and the following identities hold
\begin{eqnarray}
 \drvepsE &=&   \PRJ{\TG{\RDC}{\tsig}}   \drv{\tstrr}, \; \label{str_tangent0} \\
 \drvepsP  &=& \dps{  \PRJ{\NC{\RDC}{\tsig}}  \drv{\tstrr}},\label{str_normal} \\
 \drv{\tsig} &=&   \PRJ{\TG{\RDC}{\tsig}}\HOOK   (\drv{\tstrr}). \;  \label{const_tangent0} 
\end{eqnarray}
\end{enumerate}
When these statements are true, we also  have
\begin{eqnarray}
 \drv{\tsig}   &=&  \HOOK   (\drv{\tstrr}) -  2\mu  \PRJ{\NC{\RDC}{\tsig}}   \drv{\tstrr}, \label{projHookEpsP1} \\
  \PRJ{\NC{\RDC}{\tsig}}\HOOK   (\drv{\tstrr}) &=& 2\mu  \PRJ{\NC{\RDC}{\tsig}}   \drv{\tstrr}, \label{projHookEpsP2} \\
  \drvepsE:\drvepsP &=&  0, \label{2.14}\\
   \| \drv{\tstrr}\|^2 &=& \| \drvepsE \|^2+  \| \drvepsP \|^2. \label{2.15}
\end{eqnarray}
\end{theorem}
Before giving the proof of this theorem, let us make some comments.  An important point that emerges from this theorem is the following constitutive law
relating the stress rate and the strain rate
\begin{equation}\label{const_main_law}
 \drv{\tsig}  =  \PRJ{\TG{\RDC}{\tsig}}\HOOK   (\drv{\tstrr}).
\end{equation}
It can be observed that this law unifies the elastic and plastic regimes. 
Indeed, in the elastic regime, $\tsig$ is inside $\RDC$ and $ \PRJ{\TG{\RDC}{\tsig}} = \IDT$ 
 and we fall back onto the equations of linear elasticity. If on the other hand  $\tsig$ is on 
 the boundary of the yield domain, then  ${\TG{\RDC}{\tsig}} \ne \symt$. In this case, 
 it becomes important to give an explicit expression for the projector 
  $\PRJ{\TG{\RDC}{\tsig}}$. This will be done in Section \ref{consti_laws}
  for a yield domain defined by one or more yield functions.  
   The examples of the Von Mises criterion and the Tresca criterion will be particularly detailed in section \ref{VM-T-cr}.  \\
$\;$\\
Let us briefly describe the impact of the characterization in \eqref{const_main_law} on the 
equations of motion governing the material deformation in 
 unsteady and quasi-static regimes. The evolution equations in \eqref{evomotion_eq} and \eqref{const_main_law}
  can be gathered into the following system.
\begin{equation}\label{new_evosyst}
\frac{\pt}{\pt t}
\left(
\begin{array}{c}
\rho v \\
\tsig 
\end{array}
\right) +
\left(
\begin{array}{c}
- \div \tsig \\
  \PRJ{\TG{\RDC}{\tsig}}\HOOK (\drv{\tstrr})
\end{array}
\right) = 
\left(
\begin{array}{c}
\force \\
0
\end{array}
\right). 
\end{equation}
According to Lemma \ref{moreau_decompo} below, we have
 $\PRJ{\TG{\RDC}{\tsig}} (\ttau)+ \PRJ{\NC{\RDC}{\tsig}} (\ttau)   = \ttau$
  for all $\ttau \in \symt$.
 Combining this with \eqref{projHookEpsP2} gives 
$$
 \PRJ{\TG{\RDC}{\tsig}}\HOOK (\drv{\tstrr})  =    \HOOK (\drv{\tstrr}) -  \PRJ{\NC{\RDC}{\tsig}}  \HOOK (\drv{\tstrr}) =    \HOOK (\drv{\tstrr})  -   2\mu     \PRJ{\NC{\RDC}{\tsig}} (\drv{\tstrr}) . 
$$
Thus, 
\begin{equation}
\frac{\pt}{\pt t}
\left(
\begin{array}{c}
\rho v \\
\tsig 
\end{array}
\right) -
\left(
\begin{array}{c}
\div \tsig \\
\MHOOK (\tsig, \drv{\tstrr})
\end{array}
\right) = 
\left(
\begin{array}{c}
\force \\
0
\end{array}
\right). 
\end{equation}
with
\begin{equation}
\MHOOK(\tsig,  \ttau) = \HOOK (\ttau) -   2\mu     \PRJ{\NC{\RDC}{\tsig}}(\ttau). 
\end{equation}
In other words, from a mathematical point of view,
 the elasto-plastic incremental model of the material is obtained from
  the incremental elasticity model by replacing the linear elasticity operator 
  $\HOOK$  by the non-linear operator $\MHOOK(\tsig, .)$. Moreover,
   unlike  $\HOOK$, this non-linear operator $\MHOOK(\tsig, .)$  obviously depends 
   on the current state of stress $\tsig$.  

Time differentiation  of the first component of the system above leads to an  evolution equation involving only the velocity filed.
\begin{equation} \label{Ewaves}
\frac{\pt^2 v }{\pt t^2}  -  \div \MHOOK(\tsig, \drv{\tstrr})   = \frac{\pt h}{\pt t}. 
\end{equation}
The equations in (\ref{Ewaves})  expand the
usual elastic wave equations that are valid only within the elastic regimes. We note however that  $\MHOOK(\tsig, .)$ depends directly on $\tsig$ and in the general case, the system must be complemented with \eqref{const_main_law}.\\

For practical applications, it only remains to express the operator $ \MHOOK$ in terms of its arguments. As it will elucidated below, in the case of the Von Mises criterion \eqref{VonMisesDom}, for example, Formula \eqref{onefunc_strSIGCM} gives
\begin{equation} \label{onefunc_strSIGCM2}
\MHOOK(\tsig, \drv{\tstrr})   = \dps{ \lambda   \trc(\drv{\tstrr} ) \id+ 2\mu  \drv{\tstrr} 
                   -    \frac{\mu}{k^2} \max(0,    \drv{\tstrr}:  \dvsig)    \chi(\|\dvsig\|^2 - 2 k^2) \dvsig, }  
\end{equation}
where $\chi \; :\;  \R \to \R$ is the Heaviside type function defined by 
$\chi(t) = 1$ if $ t \geq 0$ and $\chi(t) = 0$ else.  

We now prove the theorem. 
\begin{ourproof}{of Theorem \ref{moreau_elastoplast}}
We need the following two lemmas. The first one is due to Morreau \cite{jjmoreau}. For a straightforward proof see for example   Theorem 6.29 and Corollary 6.30 of \cite{combettes_livre}. 
\begin{lemma}\label{moreau_decompo}
Let $(\HILB, \<,\>)$ be a real Hilbert space and $K$ a closed convex cone in $\HILB$. Let $K^*$ be its polar cone defined by
$$
K^* = \{\v \in \HILB \bve \forall \w \in K, \; \< \v, \w\> \leq 0\}. 
$$
Then, for any $z \in \HILB$ we have
\begin{enumerate}
\item $ \z = \PRJ{K} z + \PRJ{K^*} z$,
\item $ \< \PRJ{K} z,  \PRJ{K^*} z\> = 0$,
\item If $z = x + y$ with  $x \in K$ and $y \in   K^*$ such that $x \ne \PRJ{K} z$ and $y  \ne \PRJ{K^*} z$, then $\<x, y\> <0$. 
\end{enumerate}
\end{lemma}
\begin{lemma}\label{decompo_Hook_lem}
For all $\ttau \in \symt$ and $\tita \in \RDC$, one has
\begin{eqnarray}
 \PRJ{\TG{\RDC}{\tita}}{\HOOK   (\ttau)} &= & \HOOK (\PRJ{\TG{\RDC}{\tita}} \ttau),  \label{decompoTHook1} \\
  \PRJ{\NC{\RDC}{\tita}}{\HOOK   (\ttau)} &= & \HOOK (\PRJ{\NC{\RDC}{\tita}} \ttau),  \label{decompoNHook1}\\
  &= & 2\mu \PRJ{\NC{\RDC}{\tita}} \ttau.  \label{decompoNHook2}
\end{eqnarray}
\end{lemma}
\begin{ourproof}{} 
Let ${\ttau}_T =  \PRJ{\TG{\RDC}{\tita}}\ttau$ and $  \; \ttau_N = \PRJ{\NC{\RDC}{\tita}} \ttau$. Then,  $\ttau = \ttau_T + \ttau_N$. In view  of \eqref{indif_conesNT2} we have
\begin{equation}
\HOOK (\ttau) = ( \lambda   \trc({\ttau_T}) \id + 2 \mu {\ttau_T} ) +  2 \mu {\ttau_N}.
\end{equation}
We observe that $ \lambda   \trc({\ttau_T}) \id + 2 \mu {\ttau_T} \in \TG{\RDC}{\tita}$ (since $\TG{\RDC}{\tita}$ is a convex cone containing $\id$), 
 that $2 \mu {\ttau_N} \in  \NC{\RDC}{\tita}$ and that $( \lambda   \trc({\ttau_T}) \id + 2 \mu {\ttau_T} ):(2 \mu {\ttau_N}) = 0$. In view
 of Lemma \ref{moreau_decompo}, we deduce that $ \lambda   \trc({\ttau_T}) \id + 2 \mu {\ttau_T} = \PRJ{\TG{\RDC}{\tita}} \HOOK{\ttau}$ and
 $2 \mu {\ttau_N} =   \PRJ{\NC{\RDC}{\tita}}{\HOOK   (\ttau)}$. This ends the proof of Lemma \ref{decompo_Hook_lem}. 
\end{ourproof} 

Back to the proof of Theorem \ref{moreau_elastoplast}. We proceed in three separate steps. 
\begin{itemize}
\item[] First, we show that (1) implies (2). Assume that \eqref{sumdecompo0},  \eqref{yield_principle},  \eqref {hookzero},  \eqref{intro_normal} and \eqref{consistancy0} are satisfied. On the one hand,  since $\tsig(t) \in \RDC$ for all
$ t \in I$,  we have for all $t\in I$ and $h > 0$ sufficiently small, such that $t+h\in I$, we have $(\tsig(t+h) -\tsig(t))/h \in \TG{\RDC}{\tsig(t)}$. Taking the limit
when $h \to 0^+$ implies that $ \drv{\tsig} (t)\in \TG{\RDC}{\tsig(t)}$.  
On the other hand, \eqref{hookzero} gives 
$$
 \drvepsE(t) = \HOOK^{-1} (\drv{\tsig}(t)). 
$$
Thus, using Hooke law in (\ref{law_hooke1}) and the fact that the tangent cone contains all real multiples of the identity  (\ref{indif_conesNT}), we deduce that  $ \drvepsE  \in  \TG{\RDC}{\tsig}$.
On the other hand, we have $ \drvepsP \in \NC{\RDC}{\tsig}$, thanks to \eqref{flowrule0}.   Furthermore,
$$
\begin{array}{rcl}
\drvepsE :  \drvepsP   &=& \dps{   \HOOK^{-1} (\drv{\tsig}): \drvepsP} \\ 
&=&\dps{  \frac{1+\nu}{E} \drv{\tsig}: \drvepsP - \frac{\nu}{E} \trc{\drv{(\tsig)}} \id: \drvepsP.}
\end{array}
$$
Using the consistency condition in \eqref{consistancy0} and the fact that $\drvepsP$ is traceless according to \eqref{indif_conesNT2} yields 
$$
\drvepsE :  \drvepsP   = 0. 
$$
It follows that $\drvepsE + \drvepsP$ is the Moreau's 
decomposition of $\drv{\teps}$ described in Lemma \eqref{moreau_decompo} 
with $K =  \TG{\RDC}{\tsig}$.  This entails \eqref{str_tangent0} and
 \eqref{str_normal}. \\
Combining  \eqref{str_tangent0}, \eqref{hookzero}  and \eqref{decompoTHook1} yields
 $$
 \drv{\tsig} =   \HOOK(\drvepsE)= \HOOK( \PRJ{\TG{\RDC}{\tsig}}     \drv{\tstrr}) =  \PRJ{\TG{\RDC}{\tsig}}    \HOOK(\drv{\tstrr}). 
 $$
 Which is the identity in \eqref{const_tangent0}. 
  \item[] Second, we show that (2) implies (1). Conversely, assume that \eqref{str_tangent0},
 \eqref{str_normal} and \eqref{const_tangent0} are true and that  \eqref{yield_principle} is satisfied at $t=0$.  We need to show that \eqref{sumdecompo0}, \eqref{yield_principle}, \eqref{intro_normal} and \eqref{consistancy0} are true. 
From 
 \eqref{const_tangent0} we have hat $ \drv{\tsig} \in {\TG{\RDC}{\tsig}}$. 
  Since $\tsig(0)\in  {\RDC}$, we deduce that ${\tsig} (t)\in {\RDC}$ for all $t\in I$, meaning that  \eqref{yield_principle} is satisfied. 
  Identity \ref{sumdecompo0} follows from  \eqref{str_tangent0}, \eqref{str_normal}  and Lemma \ref{moreau_decompo}. 
  
   Combining \eqref{str_tangent0}, \eqref{const_tangent0} and \eqref{decompoTHook1} yields \eqref{hookzero}. 
For \eqref{consistancy0} we observe that 
$$
 \drvepsP : \drv{\tsig} =  \drvepsP  : \HOOK(\drvepsE) = 2 \mu  \drvepsP  :  \drvepsE+    \lambda \trc(\drvepsE ) \drvepsP : \id = 0,
$$
thanks to  Lemma \ref{moreau_decompo} and Property  \ref{indif_conesNT}. \\

\item[] Finally, identities \eqref{projHookEpsP1} and \eqref{projHookEpsP2} follow from \eqref{str_normal} and  \eqref{decompoNHook2} while \eqref{2.14} and \eqref{2.15} result from  \eqref{str_tangent0}, \eqref{str_normal}  and Lemma \ref{moreau_decompo}, and this concludes the proof of the theorem.  
\end{itemize}
\end{ourproof}
\begin{remark}
The  hydrostatic pressure invariance property in \eqref{condition_indiff_hydro} implies that
for all $\tsig \in \RDC$,  $\ttau \in \symt$,
 $\alpha > 0$ and $\beta \in \R$, we have
\begin{equation}
 \; \PRJ{\NC{\RDC}{\tsig}}(\alpha \ttau + \beta \id) = \alpha \PRJ{\NC{\RDC}{\tsig}}{\ttau}.
\end{equation}
Indeed, for  $\tita \in \NC{\RDC}{\tsig}$
$$
\|(\alpha \ttau +\beta \id) - \tita\|^2 = \alpha^2 \| \ttau - \frac{1}{\alpha} \tita\|^2 +\beta^2 \|\id\|^2 + 2 \alpha \beta \trc(\ttau), 
$$
(thanks to  \eqref{indif_conesNT2}). Therefore, minimizing $\|(\alpha \ttau +\beta \id) - \tita\|$ with respect to $ \tita$
 is equivalent to minimizing  $ \| \ttau - \alpha^{-1} \tita\|$. 
 \end{remark}

\section{Explicit constitutive laws in the case of one or two yield functions}\label{consti_laws}  
\subsection{The main practical results}
In this section we write more explicitly the constitutive laws  \eqref{str_tangent0}, \eqref{str_normal}  and \eqref{const_tangent0}  when the yield domain is defined by functional constraints. 
This is the case in most  criteria used in practice where the yield domain is often defined by a single function. 
Nevertheless, Theorems \ref{explicit_const_law_onefunc} and \ref{explicit_const_law_twofunc} stated below
deal with domains defined by several functions but the boundary points are characterized by either exactly one function and exactly two functions, respectively, 
 thus covering most practical cases. For steamlining,  the proofs of these theorems discussed in
  in sections \ref{proj_onefunc} and \ref{proj_twofunc} where the explicit calculations are presented. As one would expect, it is essentially a matter of explicitly computing the 
  projection onto the tangent and normal cones invoked in  \eqref{str_tangent0},  \eqref{str_normal}  and \eqref{const_tangent0}. 
We also note that these results are applied in Section \ref{VM-T-cr}  to the cases of  the Von Mises and Tresca criteria for illustration. Besides, the case of the Tresca criterion has  some particularities which do not fit 
completely into the framework of theorems  \ref{explicit_const_law_onefunc} and \ref{explicit_const_law_twofunc}. 
$\;$\\

As done in \cite{hanreddy}, 
\cite{koiter1953},  and \cite{koiter1960}, for example, we consider in what follows a yield domain of the form: 
\begin{equation}\label{domain_rdc}
\RDC  = \{ \tsig \in \symt \bve f_i(\tsig) \leq k_i \mbox{ for } i = 1, \cdots, m\}, \; \; 
\end{equation}
where $f_1, \cdots, f_m$ are $m$ convex differentiable functions of class $\CC^1$ and $k_1, \cdots, k_m$ are real constants. We assume that  
there exists at least one element $\tsigs \in  \symt$ such that
\begin{equation}\label{slater_condition}
f_i(\tsigs) <k_i \mbox{ for all }   i = 1, \cdots, m. 
\end{equation}
From a mathematical point of view, this condition amounts to saying that the interior of $\RDC$ is non empty and it is usually called  Slater's condition (see, e. g., \cite{borweinlewis} or 
\cite{LemareLivre}).   \\
Now, given $\tsig \in \RDC $, define $\Sat{\tsig}$ the set of indices of saturated constraints
at $\tsig$ (see, e. g.,  \cite{borweinlewis} or 
\cite{LemareLivre}): 
$$
\Sat{\tsig} = \{ i \bve 1 \leq i \leq m \mbox{ and } f_i(\tsig)  = k_i\}.
$$
The set $\Sat{\tsig}$ is empty when $\tsig$ is strictly inside the yield domain $\RDC$. 
Otherwise,  $\Sat{\tsig}  \ne \emptyset $, for all $\tsig$  on the boundary of $\RDC$. We may observe that
$$
\forall i \in \Sat{\tsig}\cap \{1,2,\cdots,m\}, \; \tgrad f_i(\tsig) \ne 0.
$$
This a straightforward consequence of the assumption in \eqref{slater_condition} because of the convexity of the functions $f_i$. 
\\

Also, in view of Slater's condition  \eqref{slater_condition},  we know that for all $\tsig \in \RDC$ such that  $\Sat{\tsig} \ne \emptyset$
$$
\begin{array}{rcl}
\TG{\RDC}{\tsig}  &=& \{ \ttau \in \symt   \bve \forall i \in \Sat{\tsig}, \ttau :  \tgrad f_i(\tsig) \leq 0 \}, \\
\NC{\RDC}{\tsig}  &=& \dps{ \{ \sum_{i \in \Sat{\tsig} } \alpha_i \tgrad f_i(\tsig)  \bve \alpha_i \geq 0 \mbox{ for all }  i  \in \Sat{\tsig} \}.}
\end{array}
$$
When $\Sat{\tsig} = \emptyset$, $\TG{\RDC}{\tsig} = \symt$ and $\NC{\RDC}{\tsig}  = \{\zero\}$.  
 We recall that we are working under the assumption that the yield functions are  insensitive to the hydrostatic pressure, that is the convex $\RDC$ in \eqref{domain_rdc} satisfies the condition in (\ref{condition_indiff_hydro}). 
A direct consequence of this assumption is that for  all $\tsig \in \RDC$, we have
\begin{equation}\label{indiff_grads}
\forall i \in \Sat{\tsig}, \; \tgrad f_i(\tsig):\id = 0. 
\end{equation}
Assumption \eqref{condition_indiff_hydro} is in particular valid when the functions $f_i, i=1,\cdots,m$, depend only on the invariants
$J_2$ and $J_3$, that is 
$$
 f_i(\tsig) = F_i(J_2(\tsig), J_3(\tsig)) \mbox{ for } 1 \leq i \leq m,
$$
where $F_1, \cdots, F_m$ are given functions of two variables. We have the following two theorems.
\begin{theorem}\label{explicit_const_law_onefunc}
Let  $\RDC$ be a yield domain defined by $m$ smooth functions as in \eqref{domain_rdc}. Assume that 
$\RDC$ satisfies the hydrostatic pressure stability condition \eqref{condition_indiff_hydro}. Let $\tsig\in \RDC$ such that $f_1(\tsig) = k_1$ and  $f_i(\tsig) < k_i$ for  all $i \geq 2$. Then, the constitutive laws in \eqref{str_tangent0},  \eqref{str_normal}, and \eqref{const_tangent0} can be rewritten as 
\begin{eqnarray}
\drv{ \tepsP}   &=&  \dps{  \frac{\max(0,\drv{ \teps}    : \tgrad f_1(\tsig))}{\|\tgrad f_1(\tsig)\|^2}  \tgrad f_1(\tsig), }\label{onefunc_strP}\\
\drv{ \tepsE}  &=& \dps{   \drv{ \teps}    - \drv{ \tepsP},  }    \label{onefunc_strE}\\ 
 \drv{\tsig} &=& \dps{  \lambda   \trc(\drv{\tstrr} ) \id+
 2 \mu \left( \drv{\tstrr}  -  \frac{\max(0,   
\drv{\tstrr}: \tgrad f_1(\tsig)    )}{\|\tgrad f_1(\tsig)\|^2}  \tgrad f_1(\tsig)\right)}. \label{constitu_law1}
\end{eqnarray}
\end{theorem}

\begin{theorem}\label{explicit_const_law_twofunc}
Let  $\RDC$ be a yield domain defined by $m$ smooth functions \eqref{domain_rdc}. Assume that 
$\RDC$ satisfies the hydrostatic pressure stability condition \eqref{condition_indiff_hydro}. Let $\tsig\in \RDC$ such that $f_i(\tsig) = k_i$ for $i=1,2$,   and that $f_i(\tsig) < k_i$ for  all $i \geq 3$. Assume also that  $ \tgrad f_1(\tsig)$ and   $ \tgrad f_2(\tsig)$ are not collinear.  
Let
$$\begin{array}{rclcrcl}
  \alpha_i &=& \dps{  \frac{ \drv{\tstrr}  : \tgrad f_i(\tsig) }{\|\tgrad f_i(\tsig)\|}  }\;\;, i = 1, 2, 
  &&       \eta_1  &=& \dps{   \frac{ \alpha_1 - \delta  \alpha_2  }{1-\delta^2}  },    \\ 
 \eta_2  &= & \dps{  \frac{ \alpha_2 - \delta \alpha_1  }{1-\delta^2} }, 
 &\mbox{ and }&    \delta &=& \dps{  \frac{  \tgrad f_1(\tsig) : \tgrad f_2(\tsig) }{  \|\tgrad f_1 (\tsig)\|   \|\tgrad f_2 (\tsig)\|}}.    
\end{array} 
$$
Then, the constitutive laws \eqref{str_tangent0},  \eqref{str_normal},  and \eqref{const_tangent0} can be rewritten as 
\begin{eqnarray}
\drv{ \tepsP} &=&
\left\{
\begin{array}{ll}
\dps{   \eta_1   \frac{ \tgrad f_1(\tsig) }{\|\tgrad f_1(\tsig)\|}  + 
\eta_2    \frac{ \tgrad f_2(\tsig) }{\|\tgrad f_2(\tsig)\|} },   & \mbox{ if } \eta_i   \geq 0 \mbox{ for } i=1, 2, \\
\dps{  \max(\alpha_1, 0)    \frac{ \tgrad f_1(\tsig) }{\|\tgrad f_1(\tsig)\|}  } ,& \mbox{ if } \min(\eta_1, \eta_2)  < 0 \mbox{ and } \alpha_1   \geq \alpha_2, \\
\dps{   \max(\alpha_2, 0)      \frac{ \tgrad f_2(\tsig) }{\|\tgrad f_2(\tsig)\|} } ,  & \mbox{ if } \min(\eta_1, \eta_2)  < 0 \mbox{ and } \alpha_1   \leq \alpha_2,
\end{array}
\right. \;\;\;\; \\
\drv{ \tepsE}  &=& \dps{   \drv{ \teps}     - \drv{ \tepsP}},  \\
\drv{\tsig} &=& \dps{  \lambda   \trc(\drv{\tstrr} ) \id+
 2 \mu  \drv{\tstrr}  - 2\mu \drv{ \tepsP}.}
\end{eqnarray}

\end{theorem}

Theorems \ref{explicit_const_law_onefunc} 
and \ref{explicit_const_law_twofunc} are part of the main results of this contribution. 
In view of Theorem \ref{moreau_elastoplast}, in order  to establish the results in these Theorems
it  is sufficient to find the expressions of the projectors  $ \PRJ{\TG{\RDC}{\tsig}}$ and  $ \PRJ{\NC{\RDC}{\tsig}}$, for yield domains as in (\ref{domain_rdc}),  in the two simple cases when $\Sat{\tsig}$ is reduced to one and two indices, respectively. This problem technical and purely computational in nature and is almost independent of mechanical modeling: it is a more general issue 
 of characterizing the projection of  any vector on the normal and tangent cones of $C$. 
 This is the subject of the following subsections. \\
 Proofs of theorems \eqref{explicit_const_law_onefunc} and \eqref{explicit_const_law_twofunc} are given in sections \ref{proj_onefunc} and \ref{proj_twofunc} hereafter. 
\subsection{A simple lemma on projections}
For $\tsig\in \RDC$ consider the subspace $\ES(\tsig)$ defined by
$$
\ES(\tsig) = \span \{ \tgrad f_i(\tsig) \bve i \in \Sat{\tsig}\},
$$
with the convention $\ES(\tsig) =  \{\zero\}$ if  $\Sat{\tsig} = \emptyset$. Denote by $\PRJP{\tsig}$ (resp. $\PRJO{\tsig}$) the orthogonal projection on 
$\ES(\tsig)$ (resp.  on $\ES(\tsig)^\perp$, where  $\ES(\tsig)^\perp$ represents the orthogonal 
subspace of $\ES(\tsig)$). 
It is easy to check that 
$$
\NC{\RDC}{\tsig} \subset \ES(\tsig)  \mbox{ and } \ES(\tsig)^\perp \subset \TG{\RDC}{\tsig}. 
$$
We have the following technical lemma.
\begin{lemma}\label{prop_carac_proj1}
Let $\tsig \in \RDC$. Then, the following identities hold true
\begin{equation}\label{identity_proj}
\PRJ{\NC{C}{\tsig}} = \PRJ{\NC{C}{\tsig}}  \circ \PRJP{\tsig} \mbox{ and } \PRJ{\TG{C}{\tsig}} = \IDT - \PRJ{\NC{C}{\tsig}}  \circ \PRJP{\tsig}.  
\end{equation}
Here, the symbol $\circ$ denotes the composition operator of functions. 
\end{lemma}
This lemma will simplify the search for the closed form expressions of the projections on the normal and tangent cones. Indeed,  given $\ttau \in \symt$, $\PRJ{\NC{C}{\tsig}} \ttau$ can be computed by following the two steps:
\begin{itemize}
\item[i.] The first step consists of finding the orthogonal projection $ \ttau^{\ES}$ of $\ttau$ on $\ES(\tsig)$: 
$$
 \ttau^{\ES} = \sum_{i \in  \Sat{\tsig}} \beta_i \tgrad f_i(\tsig), 
$$
where $\beta_i$, $i \in \Sat{\tsig}$ are (unsigned) real numbers such that
$$
 \sum_{i \in  \Sat{\tsig}} \beta_i \tgrad f_i(\tsig) : \tgrad f_k(\tsig) = \ttau: \tgrad f_k(\tsig) \mbox{ for all } k \in  \Sat{\tsig}. 
$$
\item[ii.] The second step consists of projecting $ \ttau^{\ES}$  (which lives in the finite dimensional space $\ES$) onto $\NC{\RDC}{\tsig}$. This projection coincides with
$\PRJ{\NC{C}{\tsig}} \ttau$, the projection of $\ttau$ on $\NC{C}{\tsig}$, according to Lemma \ref{prop_carac_proj1}. 
\end{itemize}
 The projection $\PRJ{\TG{C}{\tsig}}$ is  then obtained from the second identity in \eqref{identity_proj}. 
\begin{ourproof}{of Lemma \ref{prop_carac_proj1}}
Let $\ttau \in \symt$ and set $ \ttau^{\parallelsum} =  \PRJP{\tsig}  \ttau$ and
$\ttau^{\perp} = \PRJO{\tsig}  \ttau$. Thus, $\ttau = \ttau^{\parallelsum} +  \ttau^{\perp}$. In view of Lemma \ref{moreau_elastoplast}, we can also write 
$$
 \ttau^{\parallelsum}   =  \ttau_N  + \ttau_T,
$$
where $\ttau_N  = \PRJ{\NC{C}{\tsig}} \ttau^{\parallelsum}$ and $\ttau_T  = \PRJ{\TG{C}{\tsig}} \ttau^{\parallelsum}$. Set 
$\ttau^\star_T =   \ttau_T+ \ttau^{\perp}$. Obviously, $\ttau^\star_T\in \TG{C}{\tsig}$ since $\ES(\tsig)^\perp \subset \TG{\RDC}{\tsig}$ and that $\TG{C}{\tsig}$  is a convex cone.
In addition,  $\ttau^\star_T :  \ttau_N = 0$ and $\ttau =  \ttau_N  +\ttau^\star_T$. Thus, $ \ttau_N  = \PRJ{\NC{C}{\tsig}} \ttau$
and $\ttau^\star_T =  \PRJ{\TG{C}{\tsig}} \ttau$, thanks to Lemma \ref{moreau_elastoplast}. This ends the proof. 
\end{ourproof}

\begin{remark}
The hypothesis that $\RDC$ is insensitive to hydrostatic pressure is not required in Lemma \ref{prop_carac_proj1}. 
\end{remark}
\subsection{The case of one saturated yield function: Proof of Theorem \ref{explicit_const_law_onefunc} }\label{proj_onefunc}
In this subsection we consider the case of a single saturated yield function, that of points $\tsig\in \RDC$  for which
$
\card{\Sat{\tsig}} = 1, 
$
where $\card{\Sat{\tsig}}$ denotes the cardinality of the set $\Sat{\tsig}$. Without loss of generality, we assume that
 \begin{equation}\label{one_satu}f_1(\tsig) = k_1 \text{ and }f_i(\tsig) < 0 \text{ for }2 \leq i \leq m.\end{equation} Thus,  
 $$
 \TG{\RDC}{\tsig}  = \{ \tvarphi \in \symt \bve \tvarphi:\tgrad f_1(\tsig) \leq 0\},
$$
and
$$
 \NC{\RDC}{\tsig}  = \{ p \tgrad f_1(\alpha) \bve p\geq 0\}. 
$$

\begin{proposition}\label{projec_sat1}
Let $\tsig \in \RDC$ such that \eqref{one_satu} holds true. Let $\ttau \in \symt$. Then, 
\begin{eqnarray}
\PRJ{\NC{C}{\tsig}}(\ttau)  &=& \dps{ \frac{\max(0, \ttau:\tgrad f_1(\tsig))}{\|\tgrad f_1(\tsig)\|^2} \tgrad f_1(\tsig),}    \label{form_projN_onefunc}    \\
\PRJ{\TG{C}{\tsig}}(\ttau)  &=&    \ttau - \frac{\max(0, \ttau:\tgrad f_1(\tsig))}{\|\tgrad f_1(\tsig)\|^2} \tgrad f_1(\tsig).   \label{form_projT_onefunc}   
\end{eqnarray}
\end{proposition}
We recall that as a consequence of Salter's condition in \eqref{slater_condition} and convexity, we have $\tgrad f_1(\tsig) \neq 0$ whenever $f_1(\tsig)= k_1$. We may also observe that
the condition in \eqref{condition_indiff_hydro} is not needed in Proposition \ref{projec_sat1}. 
\begin{ourproof}  ~
 Let $ \ttau \in   \symt$ and set
$$
\tita = \ttau - c_0 \tgrad f_1(\tsig) \mbox{ with  }c_0  =  \frac{\max(0, \ttau:\tgrad f_1(\tsig))}{\|\tgrad f_1(\tsig)\|^2}. 
$$
Obviously $\tita  \in   \TG{C}{\tsig}$. We observe that $c_0  \tita : \tgrad f_1(\tsig) = 0$. Let $\tkappa$ be an element of $\TG{C}{\tsig}$. Then,
$$
\begin{array}{rcl}
\|\ttau - \tkappa \|^2 - \|\ttau - \tita \|^2 &= & \|(\ttau  - \tita) + (\tita - \tkappa) \|^2 - \|\ttau - \tita \|^2 \\
 &=& 2 (\ttau  - \tita):(\tita - \tkappa)+ \|\tita -\tkappa \|^2 \\
  &=& 2c_0 \tgrad f_1(\tsig):(\tita - \tkappa)+ \|\tita -\tkappa \|^2 \\
    &=& -2c_0 \tgrad f_1(\tsig):\tkappa+ \|\tita -\tkappa \|^2 \\
        &\geq& 0.  
\end{array}
$$
 It follows that  $\tita$ minimizes $\|\ttau - \tkappa \|$ over $\TG{C}{\tsig}$. 
 The rest of the proof follows thanks to lemma \ref{moreau_decompo}, from the fact that $c_0\tgrad f_1(\tsig) \in \NC{(C}{\tsig}$ and that the observation that this quantity is perpendicular to $\tita$.  
 \end{ourproof}
\begin{remark}
A sufficient condition for the yield domain $\RDC$ to be stable to hydrostatic pressure variations (Condition \ref{condition_indiff_hydro}) is that the functions $f_1,f_2,\cdots,f_m$ depend only on  the invariants $J_2$ and $J_3$. We note here that   
 if $f_1(\tsig) = F_1(J_2(\tsig), J_3(\tsig))$ for some differentiable
function $F_1$, then, 
$$
\begin{array}{rcl}
\tgrad f_1(\tsig)  =\dps{ \frac{\pt F_1}{\pt  J_2}(\IJ_2(\tsig), \IJ_3(\tsig)) \dvsig + \frac{\pt F_1}{\pt  J_3}(\IJ_2(\tsig), \IJ_3(\tsig))  (\dvsig^2  -  \frac{2}{3} \IJ_2(\tsig)\id).}
\end{array}
$$
(thanks to formula \eqref{gradJ2J3}), which yields a straightforward formula for computing the tangent and normal cones. 
\end{remark}

The proof of theorem \ref{explicit_const_law_onefunc} follows from  Formulas \eqref{form_projN_onefunc} and  \eqref{form_projT_onefunc}. 
\subsection{The case of two saturated yield function: Proof of Theorem \ref{explicit_const_law_twofunc} }\label{proj_twofunc}
In this subsection we consider the case of two saturated yield functions, that is when
$$
\card{\Sat{\tsig}} = 2. 
$$
Let's assume that
\begin{equation} 
\label{two_satu}
f_1(\tsig) = k_1, f_2(\tsig) = k_2, \text{ and  }f_i(\tsig) < k_i \text{ for } 3 \leq i \leq m. \end{equation}
We necessarily have 
\begin{equation}\label{inde_grads}
\tgrad f_1(\tsig) \mbox{ and } \tgrad f_2(\tsig)  \mbox{ are not collinear,} 
\end{equation}
otherwise, the two functions could be combined into one and we are back in the case of a single constraint.  
We know that   
 $$
 \TG{\RDC}{\tsig}  = \{ \tvarphi \in \symt \bve \tvarphi:\tgrad f_1(\tsig) \leq 0 \mbox{ and }  \tvarphi:\tgrad f_2(\tsig) \leq 0\},
$$
$$
 \NC{\RDC}{\tsig}  = \{ \eta_1 \tgrad f_1(\alpha) + \eta_2 \tgrad f_2(\alpha) \bve \eta_1\geq 0 \mbox{ and } \eta_2\geq 0\}. 
$$
As per Lemma \ref{prop_carac_proj1},  let 
$$
\ES(\tsig) = \span \{ \tgrad f_1(\tsig), \tgrad f_2(\tsig)\}. 
$$
Obviously $\NC{\RDC}{\tsig}  \subset \ES(\tsig)$.  Consider the unit vectors
$$
\ttau_i = \frac{ 1  }{\|\tgrad f_i(\tsig)\|} \tgrad f_i(\tsig), \; i=1, 2, \; 
$$
and set
$$
\delta = \ttau_1:\ttau_2 \in (-1, 1). 
$$
The orthogonal projection on $\ES(\tsig)$ is given by
\begin{equation}
\PRJP{\tsig} \ttau = \alpha_1(\ttau) \ttau_1 + \alpha_2(\ttau)  \ttau_2,
\end{equation}
where 
$$
 \alpha_1(\ttau) = \frac{(\ttau:\ttau_1)- \delta (\ttau:\ttau_2)  }{1-\delta^2} , \;  \alpha_2(\ttau) = \frac{(\ttau:\ttau_2)- \delta (\ttau:\ttau_1)  }{1-\delta^2}.
$$
\begin{proposition}\label{proj_2grads}
Under the assumptions in \eqref{two_satu} and \eqref{inde_grads}, we have, for $\ttau\in \symt$, 
\begin{equation}
\PRJ{\NC{C}{\tsig}}(\ttau)  =
\left\{
\begin{array}{ll}
 \PRJP{\tsig} \ttau & \mbox{ if }  \alpha_1(\ttau) \geq 0 \mbox{ and } \alpha_2(\ttau)  \geq 0, \\
\max(\ttau:\ttau_1, 0) \ttau_1  & \mbox{ if } \min(\alpha_1(\ttau),\alpha_2(\ttau))<0 \mbox{ and } \ttau:\ttau_1 \geq \ttau:\ttau_2, \\
 \max(\ttau:\ttau_2, 0) \ttau_2  & \mbox{ if } \min(\alpha_1(\ttau),\alpha_2(\ttau))<0 \mbox{ and }  \ttau:\ttau_1 \leq \ttau:\ttau_2.
\end{array}
\right. 
\end{equation}
\end{proposition}
\begin{ourproof}{}
Let $\ttau \in \symt$ and set $\ttau_{0} = \PRJP{\tsig} \ttau \in \ES(\tsig)$. From Lemma \ref{prop_carac_proj1} we know that
$$
\PRJ{\NC{C}{\tsig}}\ttau = \PRJ{\NC{C}{\tsig}}(\ttau_{0}).
$$ 
We have four distinct cases: 
\begin{enumerate}
\item $\ttau_{0} \in \NC{\RDC}{\tsig}$, that is  $\alpha_i(\ttau) \geq 0$  for  $ 1\leq i \leq 2$. Then, \\
 $\PRJ{\NC{C}{\tsig}} (\ttau)  = \PRJ{\NC{C}{\tsig}}(\ttau_{0}) =  \ttau_{0}$. 
\item  $\ttau_{0} \in \TG{\RDC}{\tsig}$, that is if  $\ttau:\ttau_i =  \ttau_{0}:\ttau_i \leq 0$ for  $ 1\leq i \leq 2$. Then  $\PRJ{\TG{\RDC}{\tsig}}(\ttau_{0}) =  \ttau_{0} $ and $\PRJ{\NC{\RDC}{\tsig}}(\ttau_{0}) =  \zero $. 
\item   $\ttau_{0} \not \in \NC{\RDC}{\tsig}$,  $\ttau_{0} \not  \in \TG{\RDC}{\tsig}$ and $\ttau:\ttau_1 \geq \ttau:\ttau_2$. In this case
$$
 \alpha_1(\ttau) -  \alpha_2(\ttau)  =  \frac{(\ttau:\ttau_1)-  (\ttau:\ttau_2)  }{1-\delta} \geq 0.  
$$
Thus $\alpha_1(\ttau) \geq  \alpha_2(\ttau)$.  Necessarily $\ttau:\ttau_1 > 0$   (since  $\ttau \not  \in \TG{\RDC}{\tsig}$) and $\alpha_2(\ttau) < 0$  (since  $\ttau_{0} \not \in \NC{\RDC}{\tsig}$) . We also have
$$
\PRJP{\tsig} \ttau  - (\ttau:\ttau_1) \ttau_1 = \alpha_2(\ttau)(-\delta \ttau_1+\ttau_2).
$$
Hence, $(\PRJP{\tsig} \ttau  - (\ttau:\ttau_1) \ttau_1) :\ttau_i \leq 0$ for $i =1, 2$. Thus,  $\PRJP{\tsig} \ttau  - (\ttau:\ttau_1) \ttau_1 \in  \TG{\RDC}{\tsig}$. 
Since $(\ttau:\ttau_1) \ttau_1 \in  \NC{\RDC}{\tsig}$ and $(\ttau:\ttau_1) \ttau_1 \perp \PRJP{\tsig} \ttau  - (\ttau:\ttau_1) \ttau_1$, we deduce that
 $\PRJ{\NC{\RDC}{\tsig}}(\ttau_{0})  = (\ttau:\ttau_1) \ttau_1$ and $\PRJ{\TG{\RDC}{\tsig}}(\ttau_{0}) = \PRJP{\tsig} \ttau  - (\ttau:\ttau_1) \ttau_1$. 
\item  The case $\ttau_{0} \not \in \NC{\RDC}{\tsig}$,  $\ttau_{0} \not  \in \TG{\RDC}{\tsig}$ and $\ttau:\ttau_2 \geq \ttau:\ttau_1$ is obtained by simple symmetry.
\end{enumerate}
\end{ourproof}

The proof of Theorem \ref{explicit_const_law_twofunc} follows from
Theorem \ref{moreau_elastoplast} and  Proposition \ref{proj_2grads}.
\section{Examples: Von Mises and Tresca criteria }\label{VM-T-cr}
\subsection{Von Mises criterion}
In the case of the Von Mises criterion, the yield domain is defined as
\begin{equation}\label{VonMisesDom}
\RDC = \{ \tsig \in \symt \bve \sqrt{J_2(\tsig)} \leq k\},
\end{equation}
for some given constant $k > 0$. 
\begin{corollary}
Assume that $\RDC$ is defined by \eqref{VonMisesDom} and that 
$J_2(\tsig) =  k$. Then, the constitutive laws in \eqref{str_tangent0},  \eqref{str_normal}, and \eqref{const_tangent0} 
are reduced to 
\begin{eqnarray}
\drv{ \tepsP}  &=&  \dps{  \frac{\max(0,\drv{ \teps}  : \dvsig)}{2k^2}  \dvsig, } \label{onefunc_strPCM}\\
\drv{ \tepsE}  &=& \dps{   \drv{ \teps}   -  \frac{\max(0,\drv{\teps}  : \dvsig)}{2k^2}  \dvsig},  \label{onefunc_strECM}\\
 \drv{\tsig} &=& \dps{ \lambda   \trc(\drv{\tstrr} ) \id+ 2\mu  \drv{\tstrr}
   -    \frac{\mu}{k^2} \max(0,    \drv{\tstrr}:  \dvsig)   \dvsig. }   \label{onefunc_strSIGCM}
\end{eqnarray}
\end{corollary}

\begin{ourproof}{}
This case corresponds to a single yield function 
$$
f_1(\tsig) = J_2(\tsig) = \frac{1}{2} \|\dvsig\|^2.  
$$
Using \eqref{gradJ2J3} gives
$$
\tgrad f_1(\tsig) = \dvsig.
$$
Hence,
$$
\| \tgrad f_1(\tsig) \|^2 = \| \dvsig \|^2  = 2  J_2(\tsig) = 2k^2.
$$
Replacing in \eqref{onefunc_strP}, \eqref{onefunc_strE} and \eqref{constitu_law1} gives formula 
 \eqref{onefunc_strPCM},  \eqref{onefunc_strECM} and  \eqref{onefunc_strSIGCM}.
\end{ourproof}

\subsection{Tresca criterion}\label{Tresca}
The well known Tresca criterion (or the maximum shear stress criterion) corresponds to the
yield domain
\begin{equation}\label{tresca_domain}
\RDC = \{ \tsig \in \symt \bve \ftrs(\tsig)  \leq k\},
\end{equation}
where
\begin{equation}\label{trsc_function}
\ftrs(\tsig) = \frac{1}{2} (\lda_1(\tsig) - \lda_3(\tsig)). 
\end{equation}
This is a convex function since
\begin{equation}\label{expression_ftresca}
\ftrs(\tsig) =   \frac{1}{2} (\lda_1(\tsig) + \lda_1(-\tsig)) \mbox{ and }\lda_1(\tsig)  = \max_{\|u\| = 1} \tsig: (u \otimes u).
\end{equation}
Also, the hydrostatic pressure stability condition in \eqref{condition_indiff_hydro} is satisfied since 
$$
 \forall \tsig \in \symt , \forall \theta \in \R, \; \ftrs(\tsig + \theta \id) =  \ftrs(\tsig).
$$
However, it is easy to see that the function $\ftrs$ is not everywhere differentiable.  

\begin{remark}
Of course, the definition in \eqref{expression_ftresca} is based on the assumption
 that the eigenvalues of $\tsig$ are ordered ($\lda_1(\tsig) \geq \lda_2(\tsig) \geq \lda_3(\tsig)$). However, one can also use the definition
$$
\ftrs(\tsig) = \frac{1}{4} (|\lda_1(\tsig) - \lda_2(\tsig)|+|\lda_2(\tsig) - \lda_3(\tsig)|+|\lda_1(\tsig) - \lda_3(\tsig)|),
$$
which does not depend on the order of the eigenvalues but makes it clearer that $\ftrs$ is not be everywhere differentiable.
\end{remark}

 Consistent with the condition in \eqref{condition_indiff_hydro}, we have 
\begin{equation}
\ftrs(\tsig) = \frac{1}{2} (\lda_1(\dvsig) - \lda_3(\dvsig)).
\end{equation}
We are drawn to apply the results of Proposition \ref{projec_sat1}, corresponding to the case of a single function to express the constitutive laws for the Tresca yield criterion. However, since the function is not differentiable everywhere, we will have to treat separately the points where $\ftrs$ is not differentiable. 

  In order to use Proposition \ref{explicit_const_law_onefunc}, we would like to express the yield function in terms of the components of $\tsig$ directly. One sometimes encounters in the literature the smooth function   (see, e. g., \cite{lubliner1970}, p. 137):
$$
\begin{array}{rcl}
F(\tsig) &=& \dps{ \Pi_{ 1 \leq i < j \leq 3} ((\lda_i(\tsig)-\lda_j(\tsig))^2 - 4k^2),  }\\
&=& 4 J_2(\tsig)^3 - 27 J_3(\tsig)^2 - 36 k^2 J_2(\tsig)^2 + 96 k^4 J_2(\tsig) - 64 k^6,  
\end{array}
$$
which satisfies 
$$
\ftrs(\tsig)  \leq  k \mbox{ (resp. $\ftrs(\tsig)  =  k$)} \Longrightarrow F(\tsig) \leq 0 \mbox{ (resp. $F(\tsig)  =  0$)}.
$$
This is  a one directional implication which clearly is not an equivalence. A simple way to be convinced of this, 
is to  observe that the domain $\{\tsig \bve F(\tsig) \leq 0\}$  has a smooth boundary, 
unlike the domain  \eqref{tresca_domain} defined by the Tresca criterion 
which has corners, or from a mathematical point of view, points on the boundary 
where the normal cone degenerates because the function $\ftrs$ is not differentiable there. \\
$\;$\\

It is enough to characterize the constitutive laws
 when $\tsig$ is at the boundary of the yield domain, i.e, when 
\begin{equation}\label{satu_bord_tresca}
\ftrs(\tsig)  =  k, (k>0), 
\end{equation}
i.e, when the plastic yield limit is attained. 
This characterization goes through two steps: the computation of the normal cone,
 followed by the computation of the projections on the normal   and the tangent cones. The former 
 step requires the calculation of the sub-differential of $\ftrs$.  
 
 We have two distinct cases: 
 \begin{enumerate}
\item  The three principal stresses  are distinct; i.e, $\lda_1(\tsig) > \lda_2(\tsig) > \lda_3(\tsig)$. In this case, the function 
  $\ftrs$ is differentiable at $\tsig$ (and its sub-differential is
reduced to a singleton). The projection calculation in this case falls under  Proposition \ref{projec_sat1} and Theorem \ref{explicit_const_law_onefunc}. 
\item Two of the principal stresses are equal, i.e, $\lda_2(\tsig) = \lda_1(\tsig)\ne\lda_3(\tsig) $ or
 $\lda_2(\tsig) = \lda_3(\tsig)\ne \lda_1(\tsig)$. This case is more complex because the subdifferential is not reduced to a point. The computation of the resulting projection will also be more complex
and it is not be covered by Proposition \ref{projec_sat1} for a single function   nor by Proposition  \ref{proj_2grads} for the case of two functions.
\end{enumerate}

We note that the case where the three principal stresses are equal is obviously excluded because 
of \eqref{satu_bord_tresca}.
\begin{theorem}\label{theo_vonmises}
Assume that $\RDC$ is given by \eqref{tresca_domain}. Assume that $\tsig$ satisfies \eqref{satu_bord_tresca}
and that $\lda_1(\tsig) > \lda_2(\tsig) > \lda_3(\tsig)$. 
Let $v_1(\tsig)$ and $v_3(\tsig)$ be the unitary eigenvectors associated with $\lda_1(\tsig)$ and $\lda(\tsig)$, respectively.
Then, the rules  \eqref{str_tangent0},  \eqref{str_normal}  and \eqref{const_tangent0} can be rewritten as 
\begin{eqnarray}
\drv{ \tepsP}  &=&  \dps{   q(\drv{\teps}; \tsig)  [ v_1(\tsig) \otimes v_1(\tsig) - v_3(\tsig) \otimes v_3(\tsig)], }\\
\drv{ \tepsE}  &=& \dps{   \drv{ \teps}   - \drv{ \tepsP}, } \\
 \drv{\tsig}&= &\dps{  \lambda   \trc(\drv{\tstrr} ) \id   + 2\mu\left(\drv{\tstrr} -  q(\drv{\teps}; \tsig)  [ v_1(\tsig) \otimes v_1(\tsig) - v_3(\tsig) \otimes v_3(\tsig)] \right), }
\end{eqnarray}
with 
\begin{eqnarray}
q(\drv{\teps}; \tsig) & = & \frac{1}{2} \max\left( 0,  \drv{\teps}:( v_1(\tsig) \otimes v_1(\tsig) -  v_3(\tsig) \otimes v_3(\tsig))\right), \\
 &=&   \frac{1}{2} \max\left( 0, \overline{\drv{\teps}}:( v_1(\tsig) \otimes v_1(\tsig) -  v_3(\tsig) \otimes v_3(\tsig))\right). 
\end{eqnarray}
\end{theorem}
The proof of Theorem \ref{theo_vonmises} is based on combining Proposition \ref{projec_sat1} and the following Lemma 
\begin{lemma}\label{diff_tresca_diff}
Let $\tsig \in \symt$ such that $\lda_1(\tsig) > \lda_2(\tsig) > \lda_3(\tsig)$. 
 Let $v_1(\tsig)$ and $v_3(\tsig)$ be  the unitary eigenvectors  of $\tsig$ corresponding
to  $\lda_1(\tsig)$ and  $\lda_3(\tsig)$,  respectively. Then,
\begin{equation}
\tgrad\ftrs(\tsig)= \frac{1}{2} ( v_1(\tsig) \otimes v_1(\tsig) -  v_3(\tsig) \otimes v_3(\tsig)),
\end{equation}
and thus, $\|\tgrad f (\tsig)\|^2 = 1/2$.
\end{lemma}
\begin{ourproof}{}
Since the  $\lda_i(\dvsig)$, $1 \leq i \leq 3$,  are solutions of the  equation
$$
\lambda^3 - \IJ_2(\tsig) \lambda  - \IJ_3(\tsig)  = 0,
$$
we have
\begin{eqnarray}\label{eigenvalues_tresca}
\lda_1(\dvsig) &=& \dps{  \Lambda_0 \cos \frac{\varphi_0(\tsig)}{3}, } \nonumber \\
\lda_2(\dvsig) &=& \dps{   \Lambda_0 \cos ( \frac{2\pi-\varphi_0(\tsig)}{3}), } \nonumber \\
\lda_3(\dvsig)&=& \dps{   \Lambda_0 \cos ( \frac{2\pi+\varphi_0(\tsig)}{3}).} \nonumber 
\end{eqnarray}
where 
$$
 \Lambda_0(\tsig)  = \sqrt{ \frac{4 \IJ_2(\tsig)}{3}} , \; \varphi_0(\tsig)  = \arccos \left (\frac{3\sqrt{3} \IJ_3(\tsig)}{2\IJ_2(\tsig)^{3/2}}   \right)  \in  [0, \pi]. 
$$ 
It follows that
$$
\ftrs(\tsig) = \frac{1}{2} (\lda_1(\dvsig) -  \;\lda_3(\dvsig))  =  \sqrt{\IJ_2(\tsig)}   \sin (\thet(\tsig)) \mbox{ with } \thet(\tsig) =   \frac{\pi+\varphi_0(\tsig)}{3}. 
$$
Let $v_1(\tsig)$, $v_2(\tsig)$ and $v_3(\tsig)$ be the unitary eigenvectors of $\tsig$  (and of  $\dvsig$) corresponding to
 $\lda_1(\tsig)$, $\lda_2(\tsig)$ and $\lda_3(\tsig)$, respectively. Let $P(\tsig)$  be the orthogonal matrix whose column vectors are $v_1(\tsig)$, $v_2(\tsig)$ and  $v_3(\tsig)$. Then, $\dvsig  = P(\tsig) D(\dvsig) P(\tsig)^t $ with
$$
D(\dvsig) = \diagn{\lda_1(\dvsig)}{\lda_2(\dvsig)}{\lda_3(\dvsig)} =  \frac{2 \sqrt{J_2}}{\sqrt{3}} \diagn{\cos(\alpha)}{\cos(\beta)}{\cos(\gamma)}
$$
 and $\alpha = \thet - \pi/3$, $\beta = \pi-\thet$ and $\gamma = \pi+\thet$. We have
$$
\begin{array}{rcl}
\tgrad \ftrs(\tsig) &= & \dps{   \frac{\sin (\thet(\tsig)) }{2\sqrt{ \IJ_2(\tsig) }}   \tgrad J_2(\tsig)  + \frac{1}{3} \IJ_2(\tsig) \cos(\thet(\tsig)) \tgrad \varphi_0(\tsig),} \\
&= & \dps{ \frac{\sin (\thet(\tsig)) }{ 2 \sqrt{ \IJ_2(\tsig) }}   \tgrad J_2(\tsig) }\\
&  + &\dps{\frac{\sqrt{3}\cos(\thet(\tsig))}{2 \IJ_2(\tsig)^2  \sin(\varphi_0(\tsig))}  \left(\frac{3}{2} J_3(\tsig) \tgrad J_2(\tsig) - J_2(\tsig) \tgrad J_3(\tsig)\right). }
\end{array}
$$
Since $J_3 = 2\cos(\varphi_0) J^{3/2}_2/(3\sqrt{3})$ and  $ \varphi_0 = 3\thet - \pi$, we get
\begin{equation}\label{fonc_grad}
\tgrad \ftrs(\tsig) =   \frac{\cos (2\thet) }{2\sqrt{ \IJ_2 }\sin(3\thet)}   \tgrad J_2(\tsig)  + \frac{\sqrt{3}\cos(\thet)}{2\IJ_2  \sin(3\thet)}    \tgrad J_3(\tsig).
\end{equation}
Combining with \eqref{gradJ2J3} gives
$$
\tgrad J_3(\tsig) = \frac{2}{3} J_2 P(\tsig)  \diagn{\cos(2\alpha)}{\cos(2\beta)}{\cos(2\gamma)} P(\tsig)^t. 
$$ 
Finally,
$$
\tgrad \ftrs (\tsig)= \frac{1}{2} P(\tsig)  
\diagn{1}{0}  
{-1}    P(\tsig)^t =  \frac{1}{2}  ( v_1(\tsig) v_1(\tsig)^t - v_3(\tsig) v_3(\tsig)^t).  
$$
\end{ourproof}

\begin{remark}
$\varphi_0(\tsig)/{3}$ is called the Lode angle (see \cite{wlode}). 
\end{remark}
\begin{remark}
An alternative proof of the Lemma can be formulated  using the characterization in \eqref{lewis_caract}
for the subdifferential of $\ftrs$. 
\end{remark}
We now deal with the case of a double eigenvalue. 
  Assume that  $ \lda_2(\tsig) =  \lda_\ell(\tsig)$ for some $\ell \in \{1, 3\}$ and 
  $ \lda_2(\tsig) \ne \lda_{4-\ell}(\tsig)$. Set $ m = 4 -\ell  \in \{1, 3\}$. Let $\{v_1(\tsig) , v_2(\tsig) , v_3(\tsig)\}$
 be a corresponding orthonormal basis of eigenvectors of $\tsig$.  Thus,
\begin{equation}
\begin{array}{rcl}
\tsig &=& \sum_{k=1}^3  \lda_k(\tsig) v_k(\tsig) \otimes v_k(\tsig)\\
 &=&   \lda_2(\tsig) (v_2(\tsig) \otimes v_2(\tsig) +  v_\ell(\tsig) \otimes v_\ell(\tsig))+ \lda_m(\tsig) v_m(\tsig) \otimes v_m(\tsig)\\
\end{array}
\end{equation}
Define the subspace
\begin{equation}
G_m(\tsig) = \{\ttau \in \symt \bve \ttau v_m (\tsig) = 0\}, 
\end{equation}
and set 
\begin{equation}
\begin{array}{rcl}
W_{m, i}(\tsig) &=& v_i(\tsig) \otimes v_i(\tsig), \;  \mbox{ for }  1 \leq i \leq  3 \mbox{ and } i \ne m, \\
 W_{m, m}(\tsig)  &=& \sqrt{2}\, v_2(\tsig)\odot v_\ell(\tsig). 
\end{array}
\end{equation}
We will use the following lemma whose proof is easy but reported in Appendix \ref{A1} for completeness. 
\begin{lemma}\label{lem_base_G}
Elements of $G_m(\tsig)$ commute with $\tsig$ and $\{W_{m, 1}(\tsig), W_{m, 2}(\tsig), W_{m, 3}(\tsig)\}$ is an orthonormal basis of $G_m(\tsig)$.
\end{lemma}
In what follows $\prj{G_m(\tsig)}$  denotes the orthogonal projection on  $G_m(\tsig)$. Obviously, for all   $\ttau \in \symt$
\begin{equation}
 \prj{G_m(\tsig)} \ttau = \sum_{i=1}^3 (\ttau : W_{m, i}(\tsig)) W_{m, i}(\tsig).
\end{equation}
We set
\begin{equation}
\SSD_m  (\tsig; \ttau) =  \prj{G_m(\tsig)} \dvtau - \lda_m(\prj{G_m(\tsig)}  \dvtau) \id. 
\end{equation}
We  observe that for all $\ttau \in \symt$
\begin{itemize}
\item   $\SSD_3  (\tsig; \ttau)$  (resp.  $\SSD_1  (\tsig; \ttau)$) is symmetric semidefinite positive (resp. semidefinite negative) (since $\lda_1(\prj{G_3(\tsig)} \dvtau) \geq \lda_2(\prj{G_3(\tsig)} \dvtau) \geq \lda_3(\prj{G_3(\tsig)} \dvtau)$).
\item  $\prj{G_m(\tsig)} \ttau$ and  $\SSD_m  (\tsig; \ttau)$ commute with $\tsig$.
\item $0$ is an eigenvalue of  $\prj{G_m(\tsig)} \dvtau$ (it is an eigenvalue of  all elements of $G_m(\tsig)$),
\end{itemize}
For $\prj{G_m(\tsig)} \dvtau \ne 0$, we set
\begin{equation}
\rho_m  (\tsig; \ttau)= \max(\frac{1}{4} + (-1)^{(3-m)/2} \frac{3}{4}\frac{\sum_{k=1}^3 \lda_k(\prj{G_m(\tsig)} \dvtau)}{\sum_{k=1}^3 |\lda_k(\prj{G_m(\tsig)} \dvtau)|}, 0) \in [0, 1],
\end{equation}
and, by convention, we set $\rho_m  (\tsig; \ttau)= 0$ when  $\prj{G_m(\tsig)} \dvtau = 0$.
\begin{theorem}\label{theo_tr_casegal}
Assume that $\RDC$ is given by \eqref{tresca_domain},  that $\tsig$ satisfies \eqref{satu_bord_tresca} and that  $ \lda_2(\tsig) =  \lda_\ell(\tsig)$ for some $\ell \in \{1, 3\}$ and  $ \lda_2(\tsig) \ne  \lda_{m}(\tsig)$ with $m = 4-\ell$. Then, the rules  \eqref{str_tangent0},  \eqref{str_normal}  and \eqref{const_tangent0} can be rewritten as 
\begin{eqnarray}
\drv{ \tepsP}  &=&  \dps{   \rho_m(\drv{\teps}; \tsig) [\SSD_m  (\tsig;  \drv{ \teps} )  - \trc(\SSD_m  (\tsig;  \drv{ \teps} )) v_m (\tsig) \otimes v_m(\tsig)], }\\
\drv{ \tepsE}  &=& \dps{   \drv{\teps}   - \drv{ \tepsP}, } \\
 \drv{\tsig}&= &\dps{  \lambda   \trc(\drv{\tstrr} ) \id }\\
 && \dps{   + 2\mu\left(\drv{\tstrr} -  \rho_m(\drv{\teps}; \tsig)  [\SSD_m  (\tsig;  \drv{ \teps} )  - \trc(\SSD_m  (\tsig;  \drv{ \teps} )) v_m (\tsig) \otimes v_m (\tsig)]\right). }
\end{eqnarray}
\end{theorem}
The proof of Theorem \ref{theo_tr_casegal} is mainly based on the following proposition.
\begin{proposition}\label{propo_casegal}
Assume that $f(\tsig) = k$ and that  $ \lda_1(\tsig) = \lda_2(\tsig) > \lda_3(\tsig)  $. Then,  
\begin{itemize}
\item[(a)] $ \NC{\CTRS}{\tsig} =\{ \tkappa - \trc(\tkappa ) v_3(\tsig) \otimes v_3(\tsig) \bve
 \tkappa  \in \symtp,  \tkappa  v_3(\tsig) = 0\}. $ 
 \item[(b)]  $ \NC{\CTRS}{\tsig} =\{ \tkappa - \trc(\tkappa ) v_3(\tsig) \otimes v_3(\tsig) \bve
 \tkappa  \in \symtp,  \tkappa  \tsig = \tsig \tkappa\}. $ 
\item[(c)]  $\PRJ{\NC{\CTRS}{\tsig}}(\ttau) = \PRJ{\NC{\CTRS}{\tsig}}(\dvtau)$ for all $\ttau \in \symt$.
\item[(d)]  For all $\ttau \in \symt$ 
\begin{equation}\label{formula_proj_trs_bis}
 \PRJ{\NC{\CTRS}{\tsig}}(\ttau)  = \rho_3 (\tsig; \ttau)   [\SSD_3  (\tsig; \ttau)  - \trc(\SSD_3  (\tsig; \ttau)) v_3 (\tsig) \otimes v_3 (\tsig)].
 \end{equation}
\end{itemize}
\end{proposition}
\begin{ourproof}{of Proposition \ref{propo_casegal}}
\begin{itemize}
\item[(a)] This characterization of the normal cone can be found in \cite{boulmezaoud_khouider2}. It can also be deduced from the characterization in \eqref{lewis_caract} due to \cite{lewis96} (Theorem 8.1). 
\item[(b)]   If $ \tkappa  \in \symtp,  \tkappa  v_3(\tsig) = 0$, then $0$ is an eignenvalue of $\tkappa$.   Spectral decomposition of $\tkappa$ gives  $\tkappa = \alpha_1 w_1 \otimes w_1+ \alpha_2 w_2 \otimes w_2$
with $\alpha_i \in \R$ and $w_i \in \{ v_3(\tsig)\}^\perp =  \span\{ v_1(\tsig),  v_2(\tsig)\} = \ker(\tsig - \lda_1(\tsig)\id)$, $1 \leq 1 \leq 2$. It is then easy to check that $\tsig$ and $\tkappa$ commute. \\
Conversely, if $\tsig$ and $\tkappa$ commute then $\tkappa v_3(\tsig)$ is 
an eigenvector of $\tsig$ corresponding to the eigenvalue $\lda_3$. Thus, 
  $v_3(\tsig)$ is an eigenvector of $\tkappa$. Let $\alpha_3$ be the corresponding 
  eigenvalue. Then, $\tkappa _0  = \tkappa - \alpha_3 v_3 \otimes v_3$ is also semidefinite
  positive and satisfies  $\tkappa _0 v_3(\tsig) = 0$. 
  \item[(c)] This is a direct consequence of \eqref{indif_conesNT}.
   \item[(d)] We need to calculate the projection of any  $\ttau \in \symt$ on $\NC{\CTRS}{\tsig}$. In view of the characterization above of $\NC{\CTRS}{\tsig}$ one is lead to consider the minimization problem
 \begin{equation}\label{prj_pb_trs1}
 \min_{\tkappa \in \symtp \cap G_3(\tsig)} \| \ttau - \tkappa  + \trc(\tkappa) v_3 \otimes v_3\|^2. 
 \end{equation}
 Since $\trc( \tkappa - \trc( \tkappa) v_3  \otimes  v_3) = 0$ for all $\tkappa  \in \symt$, problem \eqref{prj_pb_trs1} can be reformulated in terms of the deviatoric of $\ttau$ as follows
  \begin{equation}\label{prj_pb_trs2}
 \min_{  \tkappa \in \symtp \cap G_3(\tsig)} \| \dvtau - \tkappa + \trc(\tkappa) v_3 \otimes v_3\|^2. 
 \end{equation}
Set 
$$
F = {\rm span}\{v_3 \otimes v_3\}, \; H = (F \oplus G_3(\tsig))^\perp,
$$
where $(F \oplus G_3(\tsig))^\perp$ denote the orthogonal complement  of $F+G_3(\tsig)$ as a subspace of $\symt$. 
We may observe that $G_3(\tsig)$ is orthogonal to $F$ and 
\begin{equation}
\symt = F \oplus^\perp  G_3(\tsig)\oplus^\perp   H. 
\end{equation}
Besides, $\id  \in F \oplus G_3(\tsig)$ since 
$$
\id = \sum_{k=1}^3 v_k   \otimes v_k = v_3  \otimes v_3 +  W_1 + W_2. 
$$
Since $H = (F \oplus G_3(\tsig))^\perp$ and $\id \in F  \oplus G_3(\tsig)$ we get 
\begin{equation}\label{tita_ortho}
 \trc(\tita) =\tita: \id = 0 \mbox{ for all } \tita \in H.
\end{equation}
Now, we can write 
$$
 \dvtau  = \ttau_F + \ttau_G + \ttau_H,
$$
with $\ttau_F, \ttau_G, \; \ttau_H$ the orthogonal projections of $ \dvtau$ on the subspaces
$F$, $G$ and $H$ respectively. We may observe that
$$
\begin{array}{rcl}
\ttau_F  &=& (\dvtau : v_3  \otimes v_3)\, v_3  \otimes v_3= a_3 \, v_3  \otimes v_3 \;\;\; \mbox{ with } a_3 = \dvtau : v_3  \otimes v_3. 
\end{array}
$$
Observing that $\trc(\dvtau) = 0$ and $\trc( \ttau_H)=0$ (thanks to \eqref{tita_ortho}), 
we that  $a_3 = \trc(\ttau_F) = - \trc( \ttau_G)$.  Problem \ref{prj_pb_trs2} becomes
$$
 \min_{\tkappa \in \symtp \cap G} \| \ttau_F + \trc(\tkappa) v_3 \otimes v_3\|^2 +  \| \ttau_G - \tkappa\|^2,
$$
or 
  \begin{equation}\label{prj_pb_trs3}
 \min_{\tkappa \in \symtp \cap G_3(\tsig)}   (\trc(\tkappa)  - \trc(\ttau_G) )^2 + \| \ttau_G - \tkappa\|^2,
 \end{equation}
Let $\{w_1, w_2, v_3\}$ be an orthonormal basis of eigenvectors of $\ttau_G$ with $\mu_1$, $\mu_2$ and $0$ the corresponding  eigenvalues  (the vector $v_3$ 
is already an eigenvector). We can write
$$
\ttau_G = \mu_1 w_1 \otimes w_1 +  \mu_2 w_2 \otimes w_2. 
$$
Since $a_3 = - \trc( \ttau_G)$ we have
\begin{equation}
a_3 = - \mu_1 - \mu_2.
\end{equation}
Any tensor $\tkappa$ of $G_3(\tsig)$ can be written in the form 
$$
\tkappa = x w_1 \otimes w_1 +  y w_2 \otimes w_2 + z w_1\odot w_2,   
$$
with $x =  \tkappa  : w_1 \otimes w_1 $,  $y =  \tkappa  : w_2 \otimes w_2$ and $z = 2 \tkappa : w_1\odot w_2  = 2 \tkappa : w_1\otimes w_2$. The corresponding matrix is semidefinite positive if and only if
$$
x \geq 0, \; y \geq 0 \mbox{ and } z^2 \leq 4 xy. 
$$
In view of these elements, we can reformulate \eqref{prj_pb_trs2} as follows
  \begin{equation}\label{prj_pb_trs4}
 \min_{x \geq 0, y \geq 0, z^2 \leq xy}   (x+y -\mu_1-\mu_2)^2 +  (\mu_1-x)^2 + (\mu_2-y)^2 + z^2/4. 
 \end{equation}
We have obviously $z=0$ when the minimum is reached and the problem becomes 
  \begin{equation}\label{prj_pb_trs5}
 \min_{x \geq 0, y \geq 0}   (x+y -\mu_1-\mu_2)^2 +  (x-\mu_1)^2 + (y-\mu_2)^2. 
 \end{equation}
This quadratic optimization problem can be solved by writing usual 
Karush-Kuhn-Tucker (KKT) conditions. The minimizer
 is $\tkappa_0 = x_0 w_1 \otimes w_1 +  y_0 w_2 \otimes w_2+ z_0 w_1\odot w_2$ with 
\begin{equation}\label{formule_x0y0z0}
(x_0, y_0, z_0)  = 
\left\{
\begin{array}{ll}
(\mu_1, \mu_2, 0) & \mbox{ if } \mu_1 \geq 0  \mbox{ and }  \;\mu_2 \geq 0, \\
(0, 0, 0) & \mbox{ if } \;\mu_1 + \mu_2/2  \leq 0 \mbox{ and } \;\mu_2 + \mu_1/2 \leq 0,\\
(\mu_1 +\mu_2/2, 0, 0)  & \mbox{ if } \mu_2 < 0 \mbox{ and } \;\mu_1 + \mu_2/2 > 0, \\
(0, \mu_2 + \mu_1/2, 0)  & \mbox{ if } \mu_1 < 0 \mbox{ and } \;\mu_2 + \mu_1/2 > 0.  \\
\end{array}
\right. 
\end{equation}
We would now like to write $\tkappa_0 $ as 
$$
\tkappa_0 =   \alpha  \ttau_G + \beta (w_1 \otimes w_1 + w_2  \otimes w_2 ) = \alpha  \ttau_G + \beta (\id -  v_3  \otimes v_3).
$$
This is possible if and only if 
$$
\alpha \mu_1 + \beta =  x_0, \; \alpha \mu_2 + \beta =  y_0.
$$
Hence, if $\mu_2\neq \mu_1$ this system has one and only one solution
\begin{equation}
\alpha = \frac{y_0-x_0}{\mu_2-\mu_1},  \; \beta =  \frac{\mu_2 x_0 - \mu_1 y_0}{\mu_2-\mu_1}. 
\end{equation}
When $\mu_2=\mu_1$ we know that $x_0 = y_0$ (thanks to \eqref{formule_x0y0z0}) and we can take 
\begin{equation}
\alpha = 1 \mbox{ and } \beta = -\min(\mu_1, 0). 
\end{equation}
Using \eqref{formule_x0y0z0} we can prove that if $(\mu_1, \mu_2) \neq 0$ then
\begin{equation}
 \alpha = \max(\frac{1}{4} + \frac{3}{4}\frac{\mu_1+\mu_2}{|\mu_1|+|\mu_2|}, 0),  \beta =  -\min(\mu_1, \mu_2, 0) \alpha. 
\end{equation}
In this case,  since $\trc(M_0) = \alpha \trc(\ttau_G)  + 3 \beta$, we can write 
$$
\begin{array}{rcl}
\tkappa_0  - \trc(\tkappa_0) v_3  \otimes v_3 &=& \dps{  \alpha  \ttau_G + \beta \id  -(\alpha \trc(\ttau_G)  + 3 \beta)v_3  \otimes v_3}, \\
&=& \dps{  \alpha [ \ttau_G -  \trc(\ttau_G) v_3  \otimes v_3]     + \beta (\id  - 3 v_3  \otimes v_3)},\\
&=& \dps{  \alpha [ \ttau^\star_G -  \trc(\ttau^\star_G) v_3  \otimes v_3],}
\end{array}
$$
with $ \ttau^\star_G  =  \ttau_G -\min(\mu_1, \mu_2, 0) \id$.  This ends the proof of formula \eqref{formula_proj_trs_bis}. 
\end{itemize}
\end{ourproof}
\begin{ourproof}{of Theorem \ref{propo_casegal}}
When $ \lda_1(\tsig) = \lda_2(\tsig) > \lda_3(\tsig)$ the result is a straightforward consequence 
of Theorem \ref{moreau_elastoplast} and Proposition \ref{propo_casegal}. \\
Assume that  $ \lda_1(\tsig) > \lda_2(\tsig) = \lda_3(\tsig)$. We have $\lda_k(-\tsig) = \lda_{4-k}(\tsig)$ for $ 1 \leq k \leq 3$ and
$$
\NC{\CTRS}{-\tsig} = - \NC{\CTRS}{\tsig}.
$$
Hence  $\PRJ{\NC{\CTRS}{\tsig}}(\ttau)  = -  \PRJ{\NC{\CTRS}{-\tsig}}(-\ttau)$ and we get
$$
\PRJ{\NC{\CTRS}{\tsig}}(\ttau) = -\rho_3 (-\tsig; -\ttau)   [\SSD_3  (-\tsig; -\ttau)  - \trc(\SSD_3  (-\tsig; -\ttau)) v_3 (-\tsig) \otimes v_3 (-\tsig)].
$$
Since  $G_3(-\tsig) = G_1(\tsig)$ we deduce that
 $$
 \SSD_3  (-\tsig; -\ttau)  = - \SSD_1  (\tsig; \ttau), \;\; \rho_3 (-\tsig; -\ttau) =\rho_1 (\tsig; \ttau).  
 $$
Hence,
$$
\PRJ{\NC{\CTRS}{\tsig}}(\ttau) = \rho_1 (\tsig; \ttau)   [\SSD_1  (\tsig; \ttau)  - \trc(\SSD_1  (\tsig;\ttau)) v_1 (\tsig) \otimes v_1 (\tsig)].
$$
Combining with Theorem \ref{moreau_elastoplast}  ends the proof. 
\end{ourproof}
\begin{remark}
We have
$$
 \prj{G_3(\tsig)} \id =
  \dps{ v_1(\tsig) \otimes v_1(\tsig) +  v_2(\tsig)\otimes  v_2(\tsig)  = \id  -  v_3(\tsig) \otimes v_3(\tsig)}. 
 $$
Hence,  we get the identity
\begin{equation}
\prj{G_3(\tsig)} (\dvtau) = \prj{G_3(\tsig)} (\ttau) - \frac{\trc(\ttau)}{3}   \left(\id  -  v_3(\tsig) \otimes  v_3(\tsig)\right). 
\end{equation}
We also deduce the identities
$$
\begin{array}{rcl}
\dps{  \sum_{k=1}^3 \lda_k(\prj{G_3(\tsig)} \dvtau)  } &=& \dps{  \sum_{k=1}^3 \lda_k(\prj{G_3(\tsig)} \ttau) - \frac{2\trc(\ttau)}{3},  }\\
\dps{  \sum_{k=1}^3 |\lda_k(\prj{G_3(\tsig)} \dvtau)|  } &=&\dps{  \sum_{k=1}^3 |\lda_k(\prj{G_3(\tsig)} \ttau) - \frac{\trc(\ttau)}{3}| -  \frac{|\trc(\ttau)|}{3}. }
\end{array} 
$$
\end{remark}

 \section{Concluding remarks}
 \subsection{Main results}
 The reformulation of the perfectly elasto-plastic model described in Section \ref{sec_main_res} is based in general on two main steps: 
 \begin{itemize}
  \item[(a)]  the characterization of the normal and tangent cones at any given point of the yield surface, 
  \item[(b)]  the projection of the strain rate and its Hooke law-transform on these cones at each point of the 
  yield domain, according to Equations \ref{str_tangent0}, \ref{str_normal}, and \ref{const_tangent0}.
  \end{itemize}
  The result is an explicit evolution equation for the internal stress $\sigma$  (\ref{projHookEpsP1}) which together with the flow evolution equations form a closed system of PDEs (\ref{evol_sys}). \\
  
  In Sections \ref{consti_laws} and  \ref{VM-T-cr}, we have considered the particular case when the yield domain is described by an arbitrary (finite) number one functional constraints (convex and differentiable) and obtained the explicit expression of these projections when only one or only two of these constraints are saturated at a given boundary points.  We then looked at the most common practical examples of the Von Mises and Tresca yield criteria. As we saw, while the Von Mises examples fall into the case a single functional constraint, the example of the Tresca criterion turned out to be much more complex and required a separate treatment.  Nevertheless, but using several convex analysis and linear algebra tools, we were able to reduce it explicitly into the form in (\ref{evol_sys}). This demonstrate that our is general and can applied in principle to all imaginable yield criteria.     \\
   
   \subsection{General case of spectral yield functions}
   In the light of the Von Mises and Tresca examples considered here, it could be observed, however,
   that in most practical situations the yield functions are {\it spectral}.   A function $f \; : \;  \symt \to \R$ is said to be spectral (or weakly orthogonally invariant) if    $f(\ttau \tsig \ttau^{-1}) =  f(\tsig)$ for any $\tsig \in \symt$ and $\ttau \in \ORT{3}$ (see, e. g., \cite{lewis96} and \cite{jarre2000}).  Thus, $f$ is spectral if and only if there exists a {\it permutation invariant } function $\hf \;:\;\R^n \to \R$ such that
  \begin{equation}
  f(\tsig) = \hf(\lda_1(\tsig), \cdots, \lda_n(\tsig)). 
  \end{equation}
  ($\hf$ is said to be permutation invariant if $\hf(v_{\pi(1)}, \cdots, v_{\pi(n)}) = \hf(v)$ for any $v  \in \R^n$ and any permutation $\pi$ of $\{1, \cdots, n\}$).   Obviously, such a function $\hf$ is unique since
  \begin{equation}
\hf(v_1, \cdots, v_n) =   f(\Diag{v}). 
  \end{equation}
  For example, in the case of the Von Mises yield criterion \eqref{VonMisesDom} we have
  $$
  \hf_{M}(v_1, v_2, v_3) = \frac{1}{6} \sum_{1 \leq i < j \leq 3} (v_i - v_j)^2,
  $$
  while for the Tresca criterion \eqref{trsc_function} we have
   $$
  \hat{\ftrs}(v_1, v_2, v_3) = \frac{1}{2}  \sum_{1 \leq i < j \leq 3} |v_i -v_j|. 
  $$
Spectral functions have been the subject of much mathematical work 
which is not fully exploited in plasticity or fracture mechanics (see, e. g., \cite{lewis96}, \cite{jarre2000}, \cite{ hornBook},  \cite{dani2008} and references therein).    The purpose of this section is 
  to show that one can go much further in characterizing projections on the normal and tangent cones for spectral yield domains. 
\begin{proposition}
Assume that $\RDC$ is defined by inequalities \eqref{inequal_yield} with $f_i$ spectral for all $1 \leq i \leq m$, and that ${\rm int}{(\RDC)} \ne \emptyset$. Then, for all $\tsig$, 
 \begin{equation}\label{nrmlC_sublev}
 \begin{array}{rcl}
 \NC{\RDC}{\tsig} &=& \{  \sum_{i=1}^m  \alpha_i \ttau_i^t  \Diag{v_i}  \ttau_i  \bve  \ttau_i \in \ORT{3},  \;  \alpha_i \geq 0, \\
 &&  \; v_i \in   \pt  \hf_i (\lda(\tsig)), \alpha_i  f_i(\tsig) = 0, \   \ttau_i^t \Diag{\lda(\tsig)} \ttau_i = \tsig \}.
  \end{array}
\end{equation}
\end{proposition}
\begin{ourproof}{}
The proof is a straightforward consequence of  widely known results in convex analsyis. Set
\begin{equation}
F(\tsig) = \max_{1 \leq i \leq m} f_i(\tsig). 
\end{equation}
Obviously,  $F$ is convex, spectral and $\RDC = \{F \leq 0\}$.
If $F(\tsig)  = 0$ then  (see, e. g., \cite{LemareLivre}, Theorem 1.3.5 p 172) 
$$
 \NC{\RDC}{\tsig} = \{\mu \ttau \bve  \ttau  \in \pt F (\tsig) \mbox{ and } \mu \geq 0\}. 
$$
We also know that (see, e. g., \cite{LemareLivre}, Lemma 4.4.1)
$$
 \pt  F (\tsig) = {\rm co} \left(\cup_{i \in J(\tsig)}  \pt  f_i (\tsig)\right), \; \;\mbox{ with } J(\tsig) = \{ i \bve 1 \leq i \leq m \mbox{ and } f_i(\tsig) = 0\},
$$
(where ${\rm co(A)}$ denotes the convex hull of $A$, $A$ being a subset of $\symt$). Thus, 
$$
 \NC{\RDC}{\tsig} = \left\{\sum_{i=1}^m \mu_i \tita_i \bve \mu_i \geq 0, \; \tita_i \in   \pt  f_i (\tsig) \mbox{ and }   \mu_i f_i (\tsig)  = 0\right\}. 
$$
We then obtain \eqref{nrmlC_sublev}  by using the following characterization of $\pt  f_i$ (see \cite{lewis96}, Theorem 8.1): 
\begin{equation}\label{lewis_caract}
\begin{array}{rcl}
 \pt  f_i (\tsig)  &=& \{\ttau^t  \Diag{w}  \ttau  \bve w \in   \pt  \hf_i (\lda(\tsig)), \\
 &&          \; \ttau \in \ORT{3},\;    \ttau^t \Diag{\lda(\tsig)} \ttau = \tsig \}.
 \end{array}
\end{equation}
\end{ourproof}
\subsection{Outlook}
The formulation we proposed here clearly reveals the nature of the nonlinearity of perfect elasto-plasticity. Through the expressions in \eqref{str_tangent0} to \eqref{2.15}, the rules governing the behaviour of an elasto-plastic material  are effectively reduced to the calculation of the projectors on the tangent and normal cones of the yield domain. This led us to propose the new equation of motion 
\begin{equation}
\frac{\pt^2 v }{\pt t^2}  -  \div \MHOOK(\tsig, \drv{\tstrr} (v))  = \frac{\pt h}{\pt t}. 
\end{equation}
The study of this equation of elasto-plastic waves remains to be done.  
The authors also plan to extend the approach proposed here to
elasto-plastic deformations with hardening. This is the subject of a paper in preparation.
Moreover, a near future goal is to apply the results obtained here in the modelling and study of the behaviour of sea ice dynamics, which is shown to behave like an elasto-plastic material \cite{coon1}. However, because of the apparent technical difficulty early authors who worked on this subject either assumed that the elastic deformations are ignored when the plastic regime is reached \cite{coon1} or that the sea ice is assumed to behave as a viscous plastic material \cite{hibler}.  Also, due to large difference between the horizontal and vertical (thickness) extends of sea ice,  it is effectively considered a 2D material. 
%

\section*{Acknowledgement}
This work has been conducted in 2020-2021 when T. Z. B. was a visiting professor at the University of Victoria. This visit is funded by the French Government through the Centre National de la Recherche Scientifique (CNRS)-Pacific Institute for the Mathematical Sciences (PIMS) mobility program. The research of B. K. is partially  funded by a Natural Sciences and Engineering Research Council of Canada Discovery grant.     

\appendix 
\section{Appendix: proof of Lemma \ref{lem_base_G}} \label{A1}
Obviously,  $W_{m, i} \in G_m(\tsig)$ for $1 \leq i  \leq 3$. Furthermore, one can easily show that
for $i \ne m$ and $j \ne m$ 
$$
W_{m, i} : W_{m, j} =  \delta_{i, j}, \;  W_{m, i} : W_{m, m} = 0, \; W_{m, m} : W_{m, m}  = 1.
$$
Furthermore,  let $\tkappa \in G_m(\tsig)$. Since $\tkappa v_m(\tsig)= 0$, $0$ is an eigenvalue of $\tkappa $. Let $\{u_1, u_2, v_3\}$
be an orthonormal basis of eigenvectors of $\tkappa$. We have the spectral decomposition
$$
\tkappa  = \mu_1 u_1 \otimes u_1 +  \mu_2 u_2 \otimes u_2 + 0  u_3 \otimes u_3 = \mu_1 u_1 \mu_1 u_1 \otimes u_1 +  \mu_2 u_2 \otimes u_2. 
$$ 
On the other hand, 
$$
u_i \in (\span\{v_3\})^\perp = \span\{v_1, v_2\}  \mbox{ for } i = 1, 2.
$$
Thus, for each $i \leq 2$ there exists real numbers $\alpha_i, \beta_i$ such that
$u_i = \alpha_i v_1 + \beta_i v_2$. It follows that
$$
A  =( \mu_1 \alpha_1^2  + \mu_2 \alpha_2^2 ) W_{m, 1} + ( \mu_1 \beta_1^2  + \mu_2 \beta_2^2 ) W_{m, 2}  + 
\sqrt{2} ( \mu_1 \alpha_1 \beta_1  + \mu_2 \alpha_2 \beta_2 ) W_{m, 3}.
$$
We conclude that $\{W_1, W_2, W_3\}$ is an orthonormal basis of $G_m(\tsig)$.

 \bibliographystyle{plain}

\end{document}